\documentclass[runningheads]{llncs}
\usepackage{graphicx}
\usepackage{amsmath} % assumes amsmath package installed
\usepackage{amssymb}  % assumes amsmath package installed
\usepackage{amsfonts}
\usepackage{mathrsfs}
\usepackage{dsfont}
\newtheorem{thm}{\bf Theorem}
\newtheorem{cor}{\bf Corollary}
\newtheorem{lem}{\bf Lemma}
\newtheorem{rem}{\bf Remark}
\newtheorem{ass}{\bf Assumption}
\newtheorem{prop}{\bf Proposition}

\newtheorem{deff}{\bf Definition}
\newtheorem{eg}{\bf Example}

\newcommand{\R}{\mathbb{R}}
\newcommand{\T}{\mathbb{T}}

\newcommand{\om}{\Omega}
\newcommand{\N}{\mathbb{N}}
\newcommand{\Mm}{\mathscr{M}}
\newcommand{\W}{\mathcal{W}}
\newcommand{\U}{\mathcal{U}}
\newcommand{\M}{\mathbb{X}}
\newcommand{\G}{\mathfrak{G}}
\newcommand{\ls}{\mathcal{L}}
\newcommand{\psp}{\mathfrak{P}}
\newcommand{\Q}{\mathcal{Q}}
\newcommand{\I}{\mathbb{I}}
\newcommand{\lan}{\mathfrak{L}}
\newcommand{\ff}{\mathcal{F}}
\newcommand{\ffb}{\mathbf{F}}
\newcommand{\ob}{\mathbf{\Omega}}
\newcommand{\fff}{\mathscr{F}}
\newcommand{\pp}{\mathcal{P}}

\newcommand{\nuo}{\mathcal{\nu}_0}
\newcommand{\ppc}{\mathscr{P}}
\newcommand{\ppp}{\mathbb{P}}
\newcommand{\pim}{\mathbf{P}}

\newcommand{\ee}{\mathcal{E}}
\newcommand{\xx}{\mathcal{X}}
\newcommand{\tta}{\Theta}
\newcommand{\tmo}{\tilde{\trans{\mathcal{T}}}}
\newcommand{\tnu}{\tilde{\nu}}
 \newcommand{\tdt}{\tilde{\mathcal{T}}}
\newcommand{\w}{\varpi}
\newcommand{\eps}{\varepsilon}
\newcommand{\ep}{\vartheta}

\newcommand{\rff}{\text{ref}}
\newcommand{\bw}{\mathbb{B}_W}
\newcommand{\ttil}{\tilde{\mathcal{T}}}
\newcommand{\tv}[1]{\left\|#1\right\|_{\text{TV}}}
\newcommand{\ws}[1]{\left\|#1\right\|_{\text{W}}}
\newcommand{\wh}[1]{\widehat{#1}}
\newcommand{\fl}[1]{\lfloor{#1}\rfloor}
\newcommand{\trans}[1]{\left[\![#1\right]\!]}
\newcommand{\undert}{\check{\Theta}}
\newcommand{\overt}{\hat{\Theta}}
\newcommand{\BlackBox}{\rule{1.5ex}{1.5ex}}    
\newenvironment{pf}{\par\noindent{\bf Proof\ }}{\hfill\BlackBox}
\usepackage{dsfont}

\usepackage{hyperref}
 \hypersetup{
    colorlinks=true,
    linkcolor=blue,
    filecolor=magenta,      
    urlcolor=cyan,
    pdftitle={Overleaf Example},
    pdfpagemode=FullScreen,
    }
    \usepackage[misc]{ifsym}

\begin{document}
\title{Robustly Complete Finite-State Abstractions for Verification of Stochastic Systems}
%\thanks{This work was supported in part by the NSERC of Canada and the CRC and ERA programs.}}
%
\titlerunning{Robust Abstractions for Verification of Stochastic Systems}
\author{Yiming Meng %\orcidID{0000-1111-2222-3333}  
\and
Jun Liu%\inst{1}\orcidID{1111-2222-3333-4444} \and
%Third Author\inst{3}\orcidID{2222--3333-4444-5555
}
%\author{Omitted for submission}
%
%\authorrunning{Omitted}
\authorrunning{Y.~Meng and J.~Liu}

\institute{University of Waterloo, Waterloo, ON N2L 3G1, Canada %\and
%Springer Heidelberg, Tiergartenstr. 17, 69121 Heidelberg, Germany
\email{\{yiming.meng,\;j.liu\}@uwaterloo.ca}
}

\maketitle              % typeset the header of the contribution
\begin{abstract}
In this paper, we focus on discrete-time  stochastic systems modelled by nonlinear stochastic difference equations and propose robust abstractions for verifying probabilistic linear temporal specifications. The current literature focuses on developing sound abstraction techniques for stochastic dynamics without perturbations. However, soundness thus far has only been shown for preserving the satisfaction probability of certain types of temporal-logic specification. We present constructive finite-state abstractions for verifying probabilistic satisfaction of general $\omega$-regular linear-time properties of more general nonlinear stochastic systems. Instead of imposing  stability assumptions, we analyze the probabilistic properties from the topological view of metrizable space of probability measures. Such abstractions are both sound and approximately complete. % inverifying probabilistic satisfaction of general $\omega$-regular linear-time properties. 
 That is, given a concrete discrete-time stochastic system
and an arbitrarily small $\ls_1$-perturbation of this system, there
exists a family of finite-state  Markov chains whose set of satisfaction probabilities contains that of the original system and meanwhile is contained  by that of the slightly perturbed
system.  A direct consequence 
is that, given a probabilistic linear-time specification, initializing within the winning/losing region of the abstraction system can guarantee a satisfaction/dissatisfaction for the original system. We make an interesting observation that, unlike the deterministic case, point-mass (Dirac) perturbations cannot fulfill the purpose of robust completeness.

\keywords{Verification of stochastic systems  \and Finite-state abstraction \and Robustness \and Soundness \and Completeness \and $\ls_1$-perturbation  \and Linear temporal logic \and Metrizable space of probability measures}
\end{abstract}
\section{Introduction}
Formal verification is a rigorous mathematical technique for verifying system properties using formal analysis or model checking \cite{baier2008principles}.  So far, abstraction-based formal verification for deterministic systems has gained its maturity \cite{belta2017formal}. Whilst bisimilar (equivalent) symbolic models exist for linear (control) systems \cite{kloetzer2008fully,tabuada2006linear}, sound and approximately complete finite abstractions  can be achieved via stability assumptions \cite{pola2008approximately,girard2009approximately} or robustness (in terms of Dirac perturbations) \cite{liu2017robust,li2020robustly,liu2021closing}. 

There is a recent surge of interest in studying formal verification for stochastic systems. The verification of temporal logics for discrete-state homogeneous Markov chains can be solved by existing  tools \cite{baier2008principles,parker2013verification,bustan2004verifying,dehnert2017storm}.

In terms of verification for general discrete-time continuous-state Markov systems, a common theme is to construct abstractions to approximate the probability of satisfaction in proper ways. First attempts 
\cite{ramponi2010connections,soudjani2011adaptive,summers2010verification,abate2008probabilistic,abate2011quantitative} were to relate the 
verification of trivial probabilistic computation tree logic (PCTL) formulas to the computation of corresponding
value functions. The authors \cite{tkachev2011infinite,tkachev2012regularization} developed
alternative techniques to deal with the potential error blow-up in infinite-horizon problems. The same authors  \cite{tkachev2014characterization} investigated the necessity of  absorbing sets on the uniqueness of the solutions of corresponding Bellman equations. The related PCTL verification problem can be then precisely captured by finite-horizon ones. They also proposed abstractions for verifying general bounded linear-time (LT) properties \cite{tkachev2013formula}, and extended them to infinite-horizon reach-avoid and repeated reachability problems \cite{tkachev2013formula,tkachev2017quantitative}.

Markov set-chains are also constructive to be abstractions. The authors \cite{abate2008markov} showed that the error is finite under strong assumptions on stability (ergodicity). A closely related approach is to apply interval-valued Markov chains (IMCs), a family of finite-state Markov chains with uncertain transitions, as abstractions for the continuous-state Markov systems with certain transition kernel. 
The authors  \cite{lahijanian2015formal} argued without proof that for every PCTL formula, the probability of (path) satisfaction of the IMC abstractions forms a compact interval, which contains the real probability of the original system. They further developed `O'-maximizing/minimizing algorithms based on \cite{givan2000bounded,wu2008reachability}  to obtain the upper/lower bound of the satisfaction probability of `next', `bounded until', and `until' properties. The algorithm provides a fundamental view of computing the bounds of probability satisfaction given IMCs. However, the intuitive reasoning for soundness seems inaccurate based on our observation (readers who are interested in the details are referred to Remark \ref{rem: sound} of this paper). Inspired by \cite{lahijanian2015formal}, the work in \cite{cauchi2019efficiency} formulated IMC abstraction for verifying bounded-LTL specifications; the work in \cite{dutreix2020specification,dutreix2020verification} constructed IMC abstractions for verifying general $\omega$-regular properties of mixed-monotone systems, and provided a novel automata-based approach in obtaining the bounds of satisfaction probability. 
In \cite[Fact 1]{dutreix2020specification,dutreix2020verification}, the authors claimed the soundness of verifying general $\omega$-regular properties using IMC abstractions, but a proof is not provided. In \cite{delimpaltadakis2021abstracting}, the authors showed that IMC can be used to provide conservative estimates of the expected values for stochastic linear control system. The authors also remarked that their result can be extended to deal with omega-regular properties based on \cite{givan2000bounded,laurenti2020formal} (where  \cite{laurenti2020formal} is only regarding safety properties), but without any proofs.  To the best knowledge of the authors, we currently lacks a general framework, as the one presented in the paper, for guaranteeing soundness of IMC abstractions for verifying $\omega$-regular properties.

Motivated by these issues, our first contribution is to  provide a formal mathematical proof of the soundness of IMC abstractions for verifying $\omega$-regular linear-time properties. We show that, for any discrete-time stochastic dynamical systems modelled by a stochastic difference equation and any linear-time property, an IMC abstraction returns a compact interval of probability of (path) satisfaction which contains the satisfaction probability of the original system. A direct consequence is that starting within the winning/losing region computed by the abstraction can guarantee a satisfaction/dissatisfaction for the original system. The second contribution of this paper is to deal with stochastic systems with extra uncertain perturbations (due to, e.g., measurement errors or modeling uncertainties). Under mild assumptions, we show that, in verifying probabilistic satisfaction of general $\omega$-regular linear-time properties, IMC abstractions that are both sound and approximately complete  are constructible for nonlinear stochastic systems. That is, given a concrete discrete-time continuous-state Markov system $\M$, 
and an arbitrarily small $\ls_1$-bounded perturbation of this system, there
exists an IMC abstraction whose set of satisfaction probability contains that of $\M$, and meanwhile is contained  by that of the slightly perturbed
system. We argue in Section \ref{sec: complete} that to make the IMC abstraction robustly complete, the perturbation is generally necessary to be $\ls_1$-bounded rather than only bounded in terms of point mass. We analyze the probabilistic properties based on the topology of metrizable space of (uncertain) probability measures, and show that the technique proves more powerful than purely discussing the value of probabilities. We also would like to clarify that the main purpose of this paper is not on providing
more efficient algorithms for computing abstractions. 
We aim to provide a theoretical foundation of IMC abstractions for verifying continuous-state stochastic systems with perturbations and hope to shed some light on designing more powerful robust verification algorithms. 
%The robust model also provides a useful technical tool in the stochastic control
%design process.  
 
The rest of the paper is organized as follows. Section \ref{sec: pre} presents some preliminaries on probability spaces and Markov 
systems. Section \ref{sec: IMC} presents the soundness of abstractions in verifying $\omega$-regular linear-time properties for discrete-time nonlinear stochastic systems.  Section \ref{sec: complete} presents the constructive robust
abstractions with soundness and  approximate completeness guarantees. We discuss the differences of robustness between deterministic and stochastic systems.  The paper is concluded in
Section \ref{sec: conclusion}.

\textbf{Notation}: We denote by $\small{\prod}$ the product of ordinary sets, spaces, or function values. Denote by $\otimes$ the product of collections of sets, or sigma algebras, or measures. The $n$-times repeated product of any kind is denoted by $(\cdot)^n$ for simplification. Denote by $\pi_j:\prod_{i=0}^\infty(\cdot)_i\rightarrow(\cdot)_j$ the projection to the $j^{\text{th}}$ component. We denote by $\mathscr{B}(\cdot)$ the Borel $\sigma$-algebra of a set. 

Let $|\cdot|$ denote the inifinity norm in $\R^n$, and let $\mathbb{B}:=\{x\in\R^n: |x|\leq 1\}$. % denote the unit closed ball of $\R^n$ in infinity norm centered at the origin.  
We denote by  $\|\cdot\|_1:=\ee|\cdot|$ the $\ls_1$-norm for $\R^n$-valued random variables, and let $ \mathbb{B}_1:=\{X: \R^n\text{-valued random variable with}\;\|X\|_1\leq 1\}$. %denote by $\mathbb{B}_1$ the unit ball in $L_1$. %We denote by $\text{Lip}(\cdot)$  the Lipschitz constant  of Lipschitz continuous functions on $(\R^n, |\cdot|)$. 
Given a matrix $M$, we denote by $M_i$ its $i^{\text{th}}$ row and by $M_{ij}$ its entry at $i^{\text{th}}$ row and $j^{\text{th}}$ column. 

Given a general state space $\xx$, we denote by $\psp(\xx)$ the space of probability measures. The space of bounded and continuous functions on $\xx$ is denoted by $C_b(\xx)$.
 For any  stochastic processes $\{X_t\}_{t\geq 0}$, we use the shorthand notation $X:=\{X_t\}_{t\geq 0}$. For any stopped process $\{X_{t\wedge\tau}\}_{t\geq 0}$, where $\tau$ is a stopping time, we use the shorthand notation $X^\tau$.

\section{Preliminaries}
\label{sec: pre}
We consider $\mathbb{N}=\{0,1,\cdots\}$ as the discrete time index set, and a general Polish (complete and separable metric) space $\xx$ as the state space. For any discrete-time $\xx^\infty$-valued stochastic process $X$, we introduce some standard concepts as follows.

\subsection{Canonical sample space}\label{sec: canon}

Given a stochastic process $X$ defined on some (most likely unknown) probability space $(\om^\dagger,\mathscr{F}^\dagger, \ppp^\dagger)$.  For $\w\in \xx^\infty=:\Omega$ and $t\in \mathbb{N}$, we define $\w_t:=\pi_t(\w)$ and the coordinate process $\mathfrak{X}_t:\xx^\infty\rightarrow\xx$ as $\mathfrak{X}_t(\w):=\w_t$ associated with $\ff:=\sigma\{\mathfrak{X}_0, \mathfrak{X}_1,\cdots\}$. Then 
$\om^\dagger \longrightarrow\xx^\infty $ ($ \omega^\dagger \longmapsto \prod_{t=0}^\infty X_t(\omega^\dagger)$)
\iffalse
\begin{equation}
\begin{split}
        \om^\dagger &\longrightarrow\xx^\infty\\
        \omega^\dagger &\longmapsto \prod_{t=0}^\infty X_t(\omega^\dagger)
\end{split}
\end{equation}\fi 
is a measurable map from $(\om^\dagger,\fff^\dagger)$ to $(\om, \ff)$. In particular, $\ff=\sigma\{\mathfrak{X}_t\in \Gamma,\;\;\Gamma\in\mathscr{B}(\xx),\;\;t\in\N\}=\mathscr{B}(\xx^\infty)=\mathscr{B}^\infty(\xx)=\sigma\{\mathcal{C}\}$, where $\mathcal{C}$ is the collection of all finite-dimensional cylinder set %(in $\otimes_{t\in\N}\xx$)
of the following form:
\begin{equation*}
    \prod_{i=1}^n\Gamma_i=\{\w: \mathfrak{X}_{t_1}(\w)\in\Gamma_1, \cdots,\mathfrak{X}_{t_n}(\w)\in\Gamma_n,\;\Gamma_i\in\mathscr{B}(\xx),t_i\in \mathbb{N},i=1,\cdots,n\}.%, \;\{t_i\}_1^n\subset\N.
\end{equation*}
The measure $\pp:=\ppp^\dagger\circ X^{-1}$ of the defined coordinate process $\mathfrak{X}$ is then uniquely determined and admits the probability law of the process $X$ on the product state space, i.e.,
\begin{equation}\label{E: law}
    \begin{split}
        \pp[\mathfrak{X}_{t_1}\in\Gamma_1, \cdots,\mathfrak{X}_{t_n}\in\Gamma_n]=\pp\left(\prod_{i=1}^n \Gamma_i\right)=\ppp^\dagger[X_{t_1}\in \Gamma_1, \cdots,X_{t_n}\in\Gamma_n].
    \end{split}
\end{equation}
for any finite-dimensional cylinder set $\prod_{i=1}^n \Gamma_i\in \ff$.
We call $(\Omega,\ff,\pp)$ the canonical space of $X$ and denote by $\ee$ the associated expectation operator.
%The $\pp$ behaves like the identity to distinguish processes on $(\Omega,\ff)$.

%The canonical space gives us the convenience to study the laws of stochastic processes as well as the statistics in the state space.  
Since we only care about the probabilistic behavior of trajectories in the state space, we prefer to work on the  canonical probability spaces $(\Omega,\ff,\pp)$ and regard events as sets of sample paths. To this end, we also do not distinguish the notation $\mathfrak{X}$ from $X$ due to their identicality in distribution, i.e., we use $X$ to denote its own coordinate process for simplicity. 

In the  context of discrete state space $\xx$, we specifically use the boldface notation $(\ob,\ffb,\pim)$ for the discrete canonical spaces of some discrete-state process. 
\begin{rem}\label{rem: prod}
We usually denote by $\nu_i$ the marginal distribution of $\pp$ at some $i\in\mathbb{N}$. We can informally write the $n$-dimensional distribution (on $n$-dimensional cylinder set) as
$\pp(\cdot)=\otimes_{i=1}^n\nu_i(\cdot)$ regardless of the dependence.  %=\prod_{i=1}^n\mu_i(\pi_i\circ\cdot)$.
%$\pp[X_{t_1}\in\Gamma_1, \cdots,X_{t_n}\in\Gamma_n]=\otimes_{i=1}^n(\prod_{i=1}^n\Gamma_i)=\prod_{i=1}^n\mu_i(\Gamma_i)$. 
\end{rem}

\subsection{Markov transition systems}
For any discrete-time  stochastic process $X$, we set $\ff_t=\sigma\{X_0,X_1,\cdots,X_t\}$ to be the natural filtration. 
\begin{deff}[Markov process]
%Given a filtered probability space $(\Omega,\ff, \ff_t, \pp)$,
A stochastic process $X$ is said to be a Markov process if each $X_t$ is $\ff_t$-adapted and, for any $\Gamma\in\mathscr{B}(\xx)$ and $t>s$, we have
\begin{equation}
    \pp[X_t\in \Gamma \;|\;\ff_s]=\pp[X_t\in \Gamma\;|\;\sigma\{X_s\}],\;\;\text{a.s.}
\end{equation}
\iffalse
Equivalently, for any $\Gamma_i\in\mathscr{B}(\xx)$, $i=0,\cdots,n$,
\begin{equation}
    \pp[X_0\in\Gamma_0,\cdots,X_n\in\Gamma_n]=\int_{\cap_{i=1}^n\{X_i\in \Gamma_i\}}\pp[X_n\in\Gamma_i\;|\;X_{n-1}]d\pp(\omega)
\end{equation}
\fi
Correspondingly, for every $t$, we define the transition probability as 
\begin{equation}
    \tta_t(x,\Gamma):=\pp[X_{t+1}\in\Gamma\;|\;X_t=x],\;\;\Gamma\in\mathscr{B}(\xx).
\end{equation}
We denote $\tta_t:=\{\tta_t(x,\Gamma):\;x\in\xx,\;\Gamma\in\mathscr{B}(\xx)\}$ as the family of transition probabilities at time $t$. Note that homogeneous Markov processes are special cases such that $\tta_t=\tta_s$ for all $t\neq s$.%For $X$ with $X_0=x$, we specifically denote by $\pp^{x}$ its law.
\end{deff}
\iffalse
\begin{rem}
Homogeneous (autonomous) Markov processes are such that $\tta_t=\tta_s$ for all $t\neq s$; the $n$-step transition can be recursively defined by $\tta^{n+1}(x,\cdot)=\int_\xx\tta(y,\cdot)\tta^n(x,dy)$ with any initial distribution $\tta^0(x,\cdot)=\delta_x$.

\iffalse
\begin{equation}
    \tta^{n+1}(x,\Gamma)=\int_\xx\tta(y,\Gamma)\tta^n(x,dy),\;\;\forall A\in\mathscr{B}(\xx),
\end{equation}
\fi
\end{rem}\fi

We are interested in 
Markov processes with discrete observations of states, which is done by assigning abstract labels over a finite set of atomic propositions. We define an abstract family of labelled Markov processes as follows.
\begin{deff}[Markov system]\label{def: Markov}
%Given a filtered probability space $(\Omega,\ff, \ff_t, \pp)$,
A  Markov system is a tuple $\M=(\xx,[\![\tta]\!],\Pi,L)  $, where
\begin{itemize}
\item $\xx=\W\cup\Delta$, where $\W$ is a bounded working space, $\Delta:=\W^c$ represents all the out-of-domain states;
\item $[\![\tta]\!]$ is a collection of transition probabilities from which  $\tta_t$ is chosen for every $t$;
    %\item $x_0\in \W$ is the initial condition;
    \item $\Pi$ is the finite set of atomic propositions;
    \item $L:\xx\rightarrow 2^\Pi$ is the (Borel-measurable) labelling function.%, i.e. for every $A\in\mathscr{2^\Pi}$, $L^{-1}(A)\in\ff$.
\end{itemize}
For $X\in\M$ with $X_0=x_0$ a.s., we denote by $\pp_X^{x_0}$ %\footnote{More precisely, it should be $\pp^{X,x_0}$ for some given $X$ and $x_0$, but it is too heavy. Since we can only compute probabilities for one process w.r.t. its own law at a time, the shorthand notation should be fine.}  
the law, and $\{\pp_X^{x_0}\}_{X\in\M}$ by its collection. Similarly, for any initial distribution $\nuo\in\psp(\xx)$, we define the law by $\pp_X^{\nuo}(\cdot)=\int_\xx\pp_X^x(\cdot)\nuo(dx)$, and denote $\{\pp_X^{\nuo}\}_{X\in\M}$ by its collection. We denote by $\{\pp_{n}^{q_0}\}_{n=0}^\infty$ (resp. $\{\pp_{n}^{\nuo}\}_{n=0}^\infty$) a sequence of $\{\pp_X^{x_0}\}_{X\in\M}$ (resp. $\{\pp_X^{\nuo}\}_{X\in\M}$). We simply use $\pp_X$ (resp. $\{\pp_X\}_{X\in\M}$) if we do not emphasize the initial condition.

%We also use shorthand notation $\pp$ if the process is furthermore certain. 
\end{deff}
For a path  $\w:=\w_0\w_1\w_2\cdots\in\om$, define by 
$L_\w:=L(\w_0)L(\w_1)L(\w_2)\cdots$ its trace.  
The space of infinite words is denoted by $$(2^\Pi)^\omega=\{A_0A_1A_2\cdots:A_i\in 2^\Pi,\;i=0,1,2\cdots\}.$$ 
A linear-time (LT) property is a subset of $(2^\Pi)^\omega$. We are only interested in LT properties $\Psi$ such that $\Psi\in\mathscr{B}((2^\Pi)^\omega)$, i.e., those are Borel-measurable. 

\begin{rem}
Note that, by \cite{tkachev2017quantitative} and \cite[Proposition 2.3]{vardi1985automatic}, any $\omega$-regular language of labelled Markov processes is measurable. %The proof relies on the properties of the canonical space with the fact that $\ff=\sigma\{\mathcal{C}\}$, as well as the connection with B\"{u}chi automation. 
It follows that, for any Markov process $X$ of the  given $\M$, the traces $L_\w $  generated by measurable labelling functions are also measurable. For each $\Psi\in\mathscr{B}((2^\Pi)^\omega)$, we have the event $L_\w^{-1}(\Psi)\in\ff$. 
\end{rem}

%A linear-time (LT) property is a subset of  measurable languages. 
A particular subclass of %linear-time 
LT properties can be  specified by linear temporal logic (LTL)\footnote{While we consider LTL due to our interest, it can be easily seen that all results of this paper in fact hold for any measurable LT property, including $\omega$-regular specifications.}. 
\iffalse
This logic consists
of propositional logic operators (e.g., true, false, negation
(:), disjunction (), conjunction () and implication (!)),
and temporal operators (e.g., next (
), always (2), eventu-
ally (3), until (U) and weak until (W)).\fi
To connect with LTL specifications, we introduce the semantics of path satisfaction as well as probabilistic satisfaction as follows.
\begin{deff}
For the syntax of LTL formulae $\Psi$ and the semantics of satisfaction of $\Psi$ on infinite words, we refer readers to \cite[Section 2.4]{liu2017robust}.

For a given labelled Markov process $X$ from $\M$ with initial distribution $\nuo$, we formulate the canonical space  $(\Omega,\ff,\pp_X^{\nuo})$. For a path $\w\in\om$, we define the path satisfaction as
$$\w\vDash \Psi\Longleftrightarrow L_\w\vDash \Psi. $$
We denote by $\{ X\vDash\Psi\}:=\{\w: \;\w\vDash\Psi\}\in\ff$ %(equivalently $\{\w: \;X(\w)\vDash\Psi\}$)
 the events of path satisfaction. Given a specified probability $\rho\in[0,1]$, we define the probabilistic satisfaction of $\Psi$ as
$$X\vDash \pp^{\nuo}_{\bowtie\rho}[\Psi]\Longleftrightarrow\pp_X^{\nuo}\{X\vDash\Psi\}\bowtie \rho, $$
where $\bowtie\in\{\leq, <,\geq,>\}$.
\iffalse
For any measurable sequence $\mathfrak{L}=A_0A_1,\cdots\in\mathscr{B}((2^\Pi)^\omega)$ and a temporal logic formula $\Psi$, we inductively define $\mathfrak{L}(t) \vDash \Psi$, meaning that $\mathfrak{L}$ satisfies $\Psi$ at time $t$, as follows:
\begin{itemize}
    \item $\mathfrak{L}(t)\vDash \textbf{true}$;
    \item $\mathfrak{L}(t)\vDash \pi\Longleftrightarrow \pi\in A_i$;
    \item $\mathfrak{L}(t)\vDash \neg\Psi\Longleftrightarrow \mathfrak{L}(t) \not\vDash \Psi$;
    \item $\mathfrak{L}(t)\vDash (\Psi_1\wedge\Psi_2)\Longleftrightarrow (\mathfrak{L}(t) \vDash \Psi_1)\wedge (\mathfrak{L}(t) \vDash \Psi_2)$;
    \item $\mathfrak{L}(t)\vDash \bigcirc \Psi\Longleftrightarrow \mathfrak{L}(t+1) \vDash \Psi$;
    \item $\mathfrak{L}(t)\vDash \Psi_1\U^{\leq k}\Psi_2\Longleftrightarrow \exists s\in [t,k]$ such that $\lan(s)\vDash \Psi_2$ and $\lan(\tau)\vDash \Psi$ for all $\tau\in[t,s)$;
    \item $\mathfrak{L}(t)\vDash \Psi_1\U\Psi_2\Longleftrightarrow \exists s\geq t$ such that $\lan(s)\vDash \Psi_2$ and $\lan(\tau)\vDash \Psi$ for all $\tau\in[t,s)$;
\end{itemize}
We say $\lan$ satisfies $\Psi$, denoted by $\lan\vDash\Psi$, if $\lan(0)\vDash\Psi$. 
\fi
\end{deff}
%\begin{rem}
 %We assume that LTL formulas have been transformed into positive normal form \cite[Chapter 5]{baier2008principles}. We further assume that all negations of atomic propositions are replaced by new atomic propositions \cite{liu2017robust}.
%\end{rem}
\subsection{Weak convergence and Prokhorov's theorem}
We consider the set of possible uncertain measures within the topological space of probability measures. The following concepts are frequently used later. 
\begin{deff}[Tightness of set of measures]\label{def: tight}
Let $\xx$ be any topological state space  and $M\subseteq\psp(\xx)$ be a set of probability measures on $\xx$. 
We say that $M$ is  tight if, for every $\eps>0$ there exists
a compact set $K\subset \xx$ such that $\mu(K)\geq 1-\eps$ for every $\mu\in M$.
\end{deff}
\begin{deff}[Weak convergence]\label{def: weak_conv}
A sequence $\{\mu_n\}_{n=0}^\infty\subseteq\psp(\xx)$ is said to converge weakly to a probability measure $\mu$, denoted by $\mu_n\Rightarrow\mu$,  if 
\begin{equation}
    \int_\xx h(x)\mu_n(dx)\rightarrow \int_\xx h(x)\mu(dx),\;\;\forall h\in C_b(\xx).
\end{equation}
We frequently use the following 
alternative condition \cite[Proposition 2.2]{da2014stochastic}:
\begin{equation}
    \mu_n(A)\rightarrow\mu(A),\;\;\forall A\in\mathscr{B}(\xx) \;\text{s.t.} \;\mu(\partial A)=0.
\end{equation}

Correspondingly, the weak equivalence of any two measures $\mu$ and $\nu$ on $\xx$ is  such that 
\begin{equation}
    \int_\xx h(x)\mu(dx)=\int_\xx h(x)\nu(dx),\;\;\forall h\in C_b(\xx).
\end{equation}
\end{deff}

\begin{rem}
Weak convergence
describes the weak topology. The meaning of the  weak topology is to extend the normal convergence in deterministic settings.  Note that $x_n\rightarrow x$ in $\xx$ is equivalent to the weak convergence of Dirac measures $\delta_{x_n}\Rightarrow\delta_{x}$. It is interesting to note that $x_n\rightarrow x$ (resp. $x=y$) in $\xx$ does not imply the strong convergence (resp. equivalence) of the associated Dirac measures. A classical counterexample is to let $x_n=1/n$ and $x=0$, and we do not have $\lim_{n\rightarrow\infty}\delta_{1/n}= \delta_0$ in the strong sense since, i.e., 
$0=\lim_{n\rightarrow\infty}\delta_{1/n}(\{0\})\neq \delta_0(\{0\})=1. $

To describe the convergence (in probability law) of more general random variables $\{X_n\}$ in $\xx$, it is equivalent to investigate the weak convergence of their associated measures $\{\mu_n\}$. It is also straightforward from Definition \ref{def: weak_conv} that weak convergence also describes the convergence of probabilistic properties  related to $\{\mu_n\}$. 
\end{rem}
\begin{thm}[Prokhorov]\label{thm: proh}
Let $\xx$ be a complete separable metric space. A family $\Lambda\subseteq\psp(\xx)$ %of probability measures on $\xx$
is relatively compact %\footnote{The compactness is w.r.t. Prokhrov metric, which is uneasy to compute. We provide stronger or equivalent metrics later in this paper.} 
if an only if it is tight. Consequently, for each sequence $\{\mu_n\}$ of tight $\Lambda$, there exists a $\mu\in\bar{\Lambda}$ and a subsequence $\{\mu_{n_k}\}$ such that $\mu_{n_k}\Rightarrow\mu$.

%A family $\Lambda$ of probability measures on $\xx$ is relatively compact if and only if for any $\eps>0$, there exists a  compact set $K_\eps\subset \xx$ such that $\mu(K_\eps)\geq 1-\eps$ for all $\mu\in\Lambda$. Consequently, for each $\{\mu_n\}\subset\Lambda$, there exists a $\mu$ on $\xx$ and a subsequence $\{\mu_{n_k}\}$ such that $\mu_{n_k}\Rightarrow\mu$. 
\end{thm}
\begin{rem}
The first part of Prokhorov's theorem provides an alternative criterion for verifying the compactness of family of measures w.r.t. the corresponding metric space using tightness. On a compact metric space $\xx$, every family of probability measures is tight.
\end{rem}

\subsection{Discrete-time continuous-state stochastic systems}\label{sec: sde}
We define Markov processes determined by the difference equation
\begin{equation}\label{E: sys}
    X_{t+1}=f(X_t)+b(X_t)w_t+ \ep\xi_t%, \;\;X_0=x_0,
\end{equation}
where the state  $X_t(\w)\in\xx\subseteq\R^n$ for all $t\in\mathbb{N}$, the stochastic inputs $\{w_t\}_{t\in\mathbb{N}}$ are i.i.d.  Gaussian random variables with covariance $ I_{k\times k}$ without loss of generality. Mappings $f:\R^n\rightarrow\R^n$ and $b:\R^n\rightarrow\R^{n\times k}$ are locally Lipschitz continuous. The memoryless perturbation $\xi_t\in\mathbb{B}_1$ are independent random variables with intensity $\ep\geq 0$ %(not $\eps$) 
and unknown distributions. 

For $\ep\neq 0$,
\eqref{E: sys} defines a family $\M$ of Markov processes $X$. A special case of \eqref{E: sys} is such that $\xi$ has
Dirac (point-mass) distributions $\{\delta_x: x\in\mathbb{B}\}$ centered at some uncertain points within a unit ball.

%can be obtained by numerical schemes of stochastic differential equations with time discretizations. Discrete-time system simulates the law of the continuous-time processes only at the observation times.

\begin{rem}
The discrete-time stochastic dynamic is usually obtained from numerical schemes of stochastic differential equations driven by Brownian motions to simulate the probability laws at the observation times.  Gaussian random variables  are naturally selected to simulate Brownian motions at discrete times. Note that in \cite{dutreix2020specification},  random variables are used with known unimodal symmetric density with an interval as support. Their choice is in favor of the mixed-monotone models to provide a more accurate approximation of transition probabilities. 
Other than the precision issue, such a choice does not bring us more of the other $\ls_1$ properties. Since we focus on formal analysis based on $\ls_1$ properties rather than  providing accurate approximation, using Gaussian randomnesses as a realization does not lose any generality.
%However, this does not bring us more of the other topological and probabilistic properties. Using Gaussian randomnesses for general systems does not lose any generality.
\end{rem}
We only care about the behaviors in the bounded working space $\W$. By defining stopping time $\tau:=\inf\{t\in\mathbb{N}: X\notin\W\}$ for each $X$, we are able to study the probability law of the corresponding stopped (killed) process $X^\tau$ for any initial condition $x_0$ (resp. $\nuo$), which coincides with $\pp_X^{x_0}$ (resp. $\pp_X^{\nuo}$) on $\W$. To avoid any complexity, we use the same notation $X$ and $\pp_X^{x_0}$ (resp. $\pp_X^{\nuo}$) to denote the stopped processes and the associated laws. Such  processes driven by \eqref{E: sys} can be written as a Markov system
\begin{equation}\label{E: sysr}
    \M=(\xx, \trans{\mathcal{T}}, %x_0, 
    \Pi,L),
\end{equation}
where for all $x\in\xx\setminus\W$, the transition probability should satisfy $\mathcal{T}(x,\Gamma)=0$ for all $\Gamma\cap \W\neq \emptyset$; $\trans{\mathcal{T}}$ is the collection of transition probabilities.  For $\xi$ having Dirac distributions, the transition $\mathcal{T}$ is of the following form: 
\begin{equation}\label{E: T}
    \mathcal{T}(x,\cdot)\in\left\{\begin{array}{lr} \{\mu\sim\mathcal{N}(f(x)+\ep\xi, \;\;b(x)b^T(x)),\;\xi\in\mathbb{B} \},\;\;\forall x\in\W, \\
\{\mu:\;\mu( \Gamma)=0,\;\;\forall \Gamma\cap \W\neq \emptyset \},\;\;\forall x\in\xx\setminus\W.
\end{array}\right.
\end{equation}

\begin{ass}\label{ass: as1}
We assume that $\textbf{in}\in L(x)$ for any $x\notin\Delta$ and $\textbf{in}\notin L(\Delta)$.
We can also include `always $(\textbf{in})$' in the specifications to observe sample paths for `inside-domain' behaviors, which is equivalent to verifying $\{\tau=\infty\}$.
\end{ass}
%\begin{rem}
%For $\ep=0$, for each $x\in\W$, the $\mathcal{T}(x,\cdot)$ is a unique measure on $\xx$.
%\end{rem}
\subsection{Robust abstractions}
We define a notion of abstraction between continuous-state and finite-state Markov systems via state-level relations and measure-level relations. 
\begin{deff}
A (binary) relation $\gamma$ from $A$ to $B$ is a subset of $A\times B$ satisfying (i) for each $a\in A$, $\gamma(a):=\{b\in B: (a,b)\in \gamma\}$; (ii) for each $b\in B$, $\gamma^{-1}(b):=\{a\in A: (a,b)\in \gamma\}$; (iii) for $A'\subseteq A$, $\gamma(A')=\cup_{a\in A'}\gamma(a)$; (iv) and for $B'\subseteq B$, $\gamma^{-1}(B')=\cup_{b\in B'}\gamma^{-1}(b)$.
\end{deff}
\begin{deff}\label{def: abs}
Given a continuous-state Markov system
$$\M=(\xx, \trans{\mathcal{T}},\Pi,L) $$
and a finite-state Markov system
$$\I=(\Q, [\![\tta]\!], \Pi, L_\I),$$
where $Q=(q_1,\cdots,q_n)^T$ and $[\![\tta]\!]$ stands for a collection of $n\times n$ stochastic matrices. 
A state-level relation $\alpha\subseteq \xx\times Q$ is said to be an abstraction from $\M$ to $\I$ if (i) for all $x\in \xx$, there exists $q\in Q$ such that $(x,q)\in\alpha$; (ii) for all $(x,q)\in\alpha$, $L_\I(q)=L(x)$.

A measure-level relation
$\gamma_\alpha\subseteq \psp(\xx)\times\psp(Q)$ is said to be an abstraction from $\M$ to $\I$ if for all $i\in \{1,2,\cdots, n\}$, all $\mathcal{T}\in\trans{\mathcal{T}}$ and all $x\in\alpha^{-1}(q_i)$, there exists $\tta\in \trans{\tta}$ such that  $(\mathcal{T}(x,\cdot),\tta_i)\in\gamma_\alpha$ and that $ \mathcal{T}(x,\alpha^{-1}(q_j)) = \tta_{ij}$  for all  $j\in \{1,2,\cdots, n\}$. 

Similarly, $\gamma_\alpha\subseteq \psp(Q)\times\psp(\xx)$ is said to be an abstraction from $\I$ to  $\M$  if for all $i\in \{1,2,\cdots, n\}$, all $\tta\in \trans{\tta}$ and all $x\in\alpha^{-1}(q_i)$, there exists $\mathcal{T}\in\trans{\mathcal{T}}$ such that  $(\tta_i,\mathcal{T}(x,\cdot))\in\gamma_\alpha$ and that $ \mathcal{T}(x,\alpha^{-1}(q_j)) = \tta_{ij}$  for all $j\in \{1,2,\cdots, n\}$. 

If such relations $\alpha$ and $\gamma_\alpha $ exist, we say that $\I$ abstracts $\M$ (resp. $\M$ abstracts $\I$) and write $\M\preceq _{\gamma_\alpha} \I$ (resp. $\I\preceq _{\gamma_\alpha} \M$).
\end{deff}
\begin{ass}\label{ass: part}
Without loss of generality, we assume that the labelling function is amenable to a rectangular partition\footnote{See e.g.  \cite[Definition 1]{dutreix2020specification}.}. In other words, a state-level abstraction can be obtained from a rectangular partition.
\end{ass}

\iffalse
$\gamma_\alpha\subseteq \trans{\mathcal{T}}\times \trans{\tta}$ is said to be an abstraction from $\trans{\mathcal{T}}$ to $\trans{\tta}$ if for all $\mathcal{T}\in\trans{\mathcal{T}}$, there exists $\tta\in \trans{\tta}$ such that $(\mathcal{T},\tta)\in\gamma_\alpha$ and $ \int_{q_j}\mathcal{T}(\alpha^{-1}(q_i),dy) = \tta_{ij},\;\;\forall  q_i,q_j\in Q $. I.e.,  $\gamma_\alpha(\trans{\mathcal{T}})\subseteq \trans{\tta}$ are related through their equivalent transitions between discrete nodes. 

Similarly, $\gamma_\alpha\subseteq \trans{\tta}\times\trans{\mathcal{T}} $ is said to be an abstraction from $\trans{\tta}$ to  $\trans{\mathcal{T}}$  if for all $\tta\in \trans{\tta}$, there exists $\mathcal{T}\in\trans{\mathcal{T}}$ such that $(\tta,\mathcal{T})\in\gamma_\alpha$ and $ \int_{q_j}\mathcal{T}(\alpha^{-1}(q_i),dy) = \tta_{ij},\;\;\forall  q_i,q_j\in Q $.\fi

\section{Soundness of Robust IMC Abstractions }\label{sec: IMC}
IMCs\footnote{We omit the definition from this paper due to the limitation of space. For a formal definition see e.g. \cite[Definition 3]{lahijanian2015formal}.} are quasi-Markov systems on a discrete state space with upper/under approximations ($\overt$/$\undert$) of the real transition matrices. To abstract the transition probabilities of continuous-state Markov systems \eqref{E: sysr}, $\overt$ and $\undert$ are obtained from over/under approximations of $\mathcal{T}$ based on the state space partition. Throughout this section, we assume that  $\overt$ and $\undert$ have been correspondingly constructed.

Given an IMC, we recast it to a finite-state Markov system  
\begin{equation}\label{E: sys2}
     \I=(\Q, [\![\tta]\!], \Pi, L_\I), %\;\;\I_0=q_0\; \text{a.s.}
\end{equation}
where
\begin{itemize}
\item $\Q$ is the  finite state-space partition with dimension $N+1$ containing $\{\Delta\}$, i.e., $Q=(q_1,q_2,\cdots,q_N,\Delta)^T$;
\item $[\![\tta]\!]$\footnote{This is a necessary step to guarantee proper probability measures in \eqref{E: meas}. Algorithms can be found in \cite{hartfiel2006markov} or \cite[Section V-A]{lahijanian2015formal}.} is a set of stochastic matrices  satisfying
\begin{equation}
   [\![\tta]\!]=\{\tta: \text{stochastic matrices with}\; \undert\leq \tta\leq \overt\;\;\text{componentwisely}\};
\end{equation}
   
    \item $\Pi, L_\I$ are as before.
    
    %, i.e. for every $A\in\mathscr{2^\Pi}$, $L^{-1}(A)\in\ff$.
 %\item $q_0\in \W$ is the initial condition.
\end{itemize}
To make $\I$ an abstraction for \eqref{E: sys2}, we need the approximation to be such that $\undert_{ij}\leq \int_{q_j}\mathcal{T}(x,dy)\leq \overt_{ij}$ for all $x\in q_i$ and $i,j=1,\cdots,N$, as well as $\tta_{N+1}=(0,0,\cdots,1)$.
We further require that the partition should respect the boundaries induced by the labeling function, i.e., for any $q\in Q$,  $$L_\I(q):=L(x),\;\forall x\in q.$$
Clearly, the above connections on the state and transition probabilities satisfy Definition \ref{def: abs}. 

The Markov system $\I$ is understood as a family of `perturbed'  Markov chain generated by the uncertain choice of $\tta$ for each $t$. The $n$-step transition matrices are derived based on $[\![\tta]\!]$ as 
\begin{equation*}
    \begin{split}
      [\![\tta^{(2)}]\!]&=\{\tta_0\tta_1:\;\;\tta_0,\tta_1\in[\![\tta]\!]\},\\
      & \cdots\\
      %& \cdots\\
      [\![\tta^{(n)}]\!]&=\{\tta_0\tta_1\cdots\tta_n:\;\;\tta_i\in[\![\tta]\!],\;i=0,1,\cdots,n\}.\\
    \end{split}
\end{equation*}
Given an initial distribution $\mu_0\in\psp(Q)$, the marginal probability measure at each $t$ forms a set
\begin{equation}\label{E: meas}
    \psp(Q)\supseteq\Mm_t^{\mu_0}:=\{\mu_t=(\tta^{(t)})^T\mu_0:\;\; \tta^{(t)}\in [\![\tta^{(t)}]\!]\}. 
\end{equation}
If we do not emphasize the initial distribution $\mu_0$, we also use $\Mm_t$ to denote the marginals for short. 

We aim to show the soundness of robust IMC abstractions in this section. The proofs in this section are completed in Appendix \ref{sec: appa}.

\subsection{Weak compactness of marginal space $\Mm_t$ of probabilities }
The following lemma is rephrased from \cite[Theorem 2]{vassiliou2021non} and shows the structure of the $\Mm_t$ for each $t\in\mathbb{N}$ and any initial distribution $\mu_0$.
\begin{lem}\label{lem: compact}
Let $\I$ be a Markov system of the form \eqref{E: sys2} that is derived from an IMC. %Suppose the initial distribution is $\mu_0$. 
Then the set $\Mm_t$ of all possible probability measures at each time $t\in\mathbb{N}$ is a convex polytope, and immediately is compact. The vertices of $\Mm_t$ are of the form
\begin{equation}
    (V_{i_t})^T\cdots (V_{i_2})^T(V_{i_1})^T\mu_0
\end{equation}
for some vertices $V_{i_j}$ of $[\![\tta]\!]$.
\end{lem}
\begin{eg}
Let $Q=(q_1,q_2,q_3)^T$ and
$\mu_0=(1,0,0)^T$. The under/over estimations of transition matrices are given as
\begin{equation*}
    \undert=\begin{bmatrix}
    \frac{1}{4} & 0 & \frac{1}{4}\\
    0 & 0 & 1\\
     0 & 1 & 0\\
    \end{bmatrix},\;\;\;\overt=\begin{bmatrix}
    \frac{3}{4} & 0 & \frac{3}{4}\\
    0 & 0 & 1\\
     0 & 1 & 0\\
    \end{bmatrix}.
\end{equation*}
Then $[\![\tta]\!]$ forms a convex set of stochastic matrices with vertices
\begin{equation*}
    V_1=\begin{bmatrix}
    \frac{1}{4} & 0 & \frac{3}{4}\\
    0 & 0 & 1\\
     0 & 1 & 0\\
    \end{bmatrix},\;\;\;V_2=\begin{bmatrix}
    \frac{3}{4} & 0 & \frac{1}{4}\\
    0 & 0 & 1\\
     0 & 1 & 0\\
    \end{bmatrix}.
\end{equation*}
Therefore, the vertices of $\Mm_1$ are 
\begin{equation*}
    \nu_1^{(1)}=(V_1)^T\mu_0=(\frac{1}{4},0,\frac{3}{4})^T,\;\;\nu_1^{(2)}=(V_2)^T\mu_0=(\frac{3}{4},0,\frac{1}{4})^T.
\end{equation*}
Hence, 
 $\Mm_1=\{\mu:\mu=\alpha\nu_1^{(1)}+(1-\alpha)\nu_1^{(2)},\;\alpha\in[0,1]\}$.
 Similarly, the vertices of $\Mm_2$ are
 \begin{equation*}
     \begin{split}
        & \nu_2^{(1)}=(V_1)^T(V_1)^T\mu_0=(\frac{1}{16},\frac{12}{16},\frac{3}{16})^T,\;\;\nu_2^{(2)}=(V_2)^T(V_1)^T\mu_0=(\frac{3}{16},\frac{12}{16},\frac{1}{16})^T,\\
        & \nu_2^{(3)}=(V_1)^T(V_2)^T\mu_0=(\frac{3}{16},\frac{4}{16},\frac{9}{16})^T,\;\;\nu_2^{(4)}=(V_2)^T(V_2)^T\mu_0=(\frac{9}{16},\frac{4}{16},\frac{3}{16})^T,\\
     \end{split}
 \end{equation*}
 and 
 \begin{equation*}
     \Mm_2=\{\mu:\mu=\alpha\beta\nu_2^{(1)}+\alpha(1-\beta)\nu_2^{(2)}+\beta(1-\alpha)\nu_2^{(3)}+(1-\alpha)(1-\beta)\nu_2^{(4)},\;\alpha,\beta\in[0,1]\}.
 \end{equation*}
 The calculation of the rest of $\Mm_t$ follows the same procedure.
\end{eg}
Now we introduce the total variation distance $\tv{\;\cdot\;}$ and see how $(\Mm_t,\tv{\;\cdot\;})$ (at each $t$) implies the weak topology.
\begin{deff}[Total variation distance]
Given two probability measures $\mu$ and $\nu$ on $\xx$, the total variation distance is defined as
\begin{equation}
    \tv{\mu-\nu}=2\sup_{\Gamma\in\mathscr{B}(\xx)}|\mu(\Gamma)-\nu(\Gamma)|.
\end{equation}
    In particular, if $\xx$ is a discrete space, 
    \begin{equation}
         \tv{\mu-\nu}=\|\mu-\nu\|_1=\sum_{q\in\xx}|\mu(q)-\nu(q)|.
    \end{equation}
\end{deff}
\iffalse
\begin{rem}
It is equivalent to use the dual representation
$$\tv{\mu-\nu}=\sup_{\substack{ \|h\|_\infty\leq 1} }\left|\int_\xx h(x)\mu(dx)-\int_\xx h(x)\nu(dx)\right|.$$
\end{rem}\fi

\begin{cor}\label{cor: marg}
Let $\I$ be a Markov system of the form \eqref{E: sys2} that is derived from an IMC. Then at each time $t\in\mathbb{N}$,  for for each $\{\mu_n\}\subseteq\Mm_t$, there exists a $\mu\in\Mm_t$ and a subsequence $\{\mu_{n_k}\}$ such that $\mu_{n_k}\Rightarrow\mu$. In addition, for each $h\in C_b(\xx)$ and $t\in\mathbb{N}$, the set $H=\{\sum_\xx h(x)\mu(x),\;\mu\in\Mm_t\}$ forms a convex and compact subset in $\R$.
\end{cor}

\begin{rem}\label{rem: tight}
The above shows that $\tv{\;\cdot\;}$ metricizes the weak topology of measures on $Q$. Note that since $Q$ is bounded and finite, any metrizable family of measures on $Q$ is compact. For example, let $Q=\{q_1,q_2\}$, and $\{(0,1)^T,(1,0)^T\}$ be a set of singular measures on $Q$. Then every sequence of the above set has a weakly convergent subsequence. However, these measures do not have a convex structure as $\Mm_t$. Hence the corresponding $H$ that is generated by $\{(0,1)^T,(1,0)^T\}$ only provides vertices in $\mathbb{Z}$.
%given that $Q$ is a bounded and finite state space. 
\end{rem}
\subsection{Weak compactness of probability laws of $\I$ on infinite horizon} 
In this subsection, we focus on the case where $I_0=q_0$ a.s. for any $q_0\in Q\setminus \{\Delta\}$. The cases for arbitrary initial distribution should be similar.
We formally denote $\Mm^{q_0}:=\{\pim_I^{q_0}\}_{I\in\I}$ by the set of probability laws of every discrete-state Markov processes $I\in\I$ with initial state $q_0\in Q$. We denote $\Mm_t^{q_0}$ by the set of marginals at $t$.
%We use $\pim$ to denote elements of $\Mm$. 
\begin{prop}\label{prop: law_weak}
For any $q_0\in Q$, every sequence $\{\pim_{n}^{q_0}\}_{n=0}^\infty$ of $\Mm^{q_0}$ has a weakly  convergent subsequence.
\end{prop}
\begin{rem}
The property is an extension of the marginal weak compactness relying on the (countable) product topology. We can also introduce proper product metrics to metricize, see e.g. \cite{sagar2021compactness}. Similar results hold under certain conditions for continuous time processes on continuous state spaces with uniform norms \cite[Lemma 82.3 and 87.3]{rogers2000diffusions}. 
\end{rem}

\begin{thm}\label{thm: imc_bound}
Let $\I$ be a Markov system of the form \eqref{E: sys2} that is derived from an IMC. %Let $q_0\in Q$ be the initial state. 
%Let $\{\iim\}$ be the set of all the discrete-state non-homogeneous Markov processes of $\I$ and $\Mm=\{\pim^\iim\}$ be the corresponding set of laws of $\I$. 
Then for any LTL formula $\Psi$, the set $S^{q_0}=\{\pim_I^{q_0}(I\vDash\Psi)\}_{I\in\I}$ is a convex and compact subset in $\R$, i.e., a compact interval. 
\end{thm}

\subsection{Soundness of IMC abstractions}
%Recall notations from \ref{sec: sde}.
\begin{prop}\label{cor: compact}
Let $\M$ be a Markov system  driven by \eqref{E: sysr}. Then every sequence $\{\pp_n^{x_0}\}_{n=0}^\infty$ of $\{\pp_X^{x_0}\}_{X\in\M}$ has a weakly  convergent subsequence. Consequently, for any LTL formula $\Psi$, the set $\{\pp_X^{x_0}(X\vDash\Psi)\}_{X\in\M}$ is a compact subset in $\R$.
\end{prop}

\begin{lem}\label{lem: connect}
Let $X\in\M$ be any Markov process  driven by \eqref{E: sysr} and $\I$ be the finite-state IMC abstraction of $\M$. Suppose the initial distribution $\nu_0$ of $X$ is such that $ \nu_0(q_0)=1$. Then, 
 there exists a unique law $\pim_I^{q_0}$ of some $I\in\I$   such that, for any LTL formula $\Psi$,
$$\pp_X^{\nu_0}(X\vDash\Psi)=\pim_I^{q_0}(I\vDash\Psi).$$
\end{lem}

\begin{thm}\label{thm: inclusion}
Assume the settings in Lemma \ref{lem: connect}. For any LTL formula $\Psi$, we have$$ \pp_X^{\nu_0}(X\vDash\Psi)\in \{\pim_I^{q_0}(I\vDash\Psi)\}_{I\in\I},$$
%where $S$, as given in Theorem \ref{thm: imc_bound}, is a convex and compact subset in $\R$.
\end{thm}
\begin{pf}
The conclusion is obtained by combining Lemma \ref{lem: connect} and Theorem \ref{thm: imc_bound}.
\end{pf}
\begin{cor}
Let $\M$, its IMC abstraction $\I$, an LTL formula $\Psi$, and a constant $\rho\in[0,1]$ be given. Suppose $I\vDash\pim^{q_0}_{\bowtie\rho}[\Psi]$ for all $I\in\I$, we have $X\vDash\pp^{\nu_0}_{\bowtie\rho}[\Psi]$ for all $X\in\M$ with $\nu_0(q_0)=1$. 
\end{cor}
\begin{rem}\label{rem: sound}
%We can simply focus on the special cases where $\nu_0$ is $\delta_{x_0}$ for any $x_0\in q_0$. 
Note that  we do not have $\pp_X^{\nu_0}\in\{\pim_I^{q_0}\}_{I\in\I}$ since each $\pim_I^{q_0}$ is a discrete measure whereas $\pp_X^{\nu_0}$ is not. They only coincide when measuring Borel subset of $\ffb$. It would be more accurate to state that $\pp_X^{\nu_0}(X\vDash\Psi)$ is a member of $\{\pim_I^{q_0}(I\vDash\Psi)\}_{I\in\I}$  rather than say ``the true distribution (the law as what we usually call) of the original system is a member of
the distribution set represented by the abstraction model'' \cite{lahijanian2015formal}. 
\end{rem}
\begin{rem}
We have seen that, in view of Lemma \ref{lem: connect}, the `post-transitional' measures are automatically related only based on the relations between  transition probabilities. We will see in the next section that such relations can be constructed to guarantee an approximate completeness of $\I$.
\end{rem}

\begin{prop}\label{prop: converge}
Let $\eps:=\max_i\|\overt_i-\undert_i\|_{\text{TV}}$. %, where $i $ is the index for rows of matrices.
Then for each LTL formula $\Psi$, as $\eps\rightarrow 0$, the length $\lambda(S^{q_0})\rightarrow 0$.
\end{prop}

\begin{rem}
By Lemma \ref{lem: connect}, for each $X\in\M$, there exists exactly one $\pim_I$ of some $I\in\I$ by which  satisfaction probability equals to that of $X$.
%IMC abstractions are nothing but to approximate the exact probability of satisfaction. 
The precision of $\overt$ and $\undert$ determines the size of $S^{q_0}$. 
Once we are able to calculate the exact law of $X$, the $S^{q_0}$ becomes a singleton by Proposition \ref{prop: converge}. For example,
 let each $w_t$ become $\delta_0$, we have each $\Mm_t$ reduced to a singleton $\{\delta_{f(x_t)}\}$ automatically. The verification problem becomes checking whether $L(f(x_t))\vDash\Psi$ given the partition $Q$. The probability of satisfaction is either $0$ or $1$. Another example would be $X_{t+1}=AX_t+Bw_t$, where $A,B$ are linear matrices. We are certain about the exact law of this system, and there is no need to introduce IMC for approximations at the beginning. IMC abstractions prove more useful when coping with systems whose marginal distributions are uncertain or not readily computable. 
\end{rem}

\section{Robust Completeness of IMC Abstractions}\label{sec: complete}
In this section, we are given a Markov system $\M_1$ driven by \eqref{E: sys} with  point-mass perturbations of strength $\ep_1\geq 0$. Based on $\M_1$, we first construct an IMC abstraction $\I$. We then show that $\I$ can be abstracted by a system  $\M_2$ with more general $\ls_1$-bounded noise of any arbitrary strength $\ep_2>\ep_1$.
%we show that robust IMC abstractions for the (perturbed) discrete-time stochastic system \eqref{E: sysr} is approximately complete. That is, for arbitrary numbers $0\leq\ep_1<\ep_2$ and the corresponding perturbed systems $\M_1$ and $\M_2$, there is a finite-state IMC $\I$ that abstracts $\M_1$ and meanwhile is abstracted by $\M_2$. 

\iffalse
\begin{rem}
We have seen that $\M$ and the abstraction $\I$ are connected through  a two-level relation, i.e., a state-level $\alpha\subseteq \xx\times Q$ and a measure(law)-level $\beta=\{(\nu,\mu)\in\ppc\times\Mm: \nu_t(\cup_{x\in\alpha^{-1}(q)}x)=\mu_t(q),\;\forall q\in Q,\forall t\}$ such that $\beta(\ppc)\subseteq\Mm$.  In view of Lemma \ref{lem: connect}, it is clear that $\beta$ is induced only by the transition probabilities requiring that $\cup_{x\in\alpha^{-1}(q_i)}\mathcal{T}(x,q_j)\subseteq \trans{\tta_{ij}}$, which is exactly the requirement of $\I$. We will see in the next section that a relation between transition probabilities of $\M$ and $\I$  can be constructed to guarantee the approximate completeness of $\I$ for probabilistic verification.
\end{rem}\fi

 Recalling the soundness analysis of IMC abstractions in Section \ref{sec: IMC}, the relation of satisfaction probability is induced by a relation between the continuous and discrete transitions. To capture the probabilistic properties of stochastic processes, reachable set of probability measures is the analogue of the reachable set in deterministic cases. We rely on a similar technique in this section to discuss how transition probabilities of different uncertain Markov systems are related. To metricize sets of Gaussian measures and to connect them with discrete measures, we prefer to use Wasserstein metric. 
\begin{deff}
Let $\mu,\nu\in\psp(\xx)$ for $(\xx,|\cdot|)$, the Wasserstein distance\footnote{This is formally termed as $1^{\text{st}}$-Wasserstein metric. We choose $1^{\text{st}}$-Wasserstein metric due to the convexity and nice property of test functions.} is  defined by
$
    \ws{\mu-\nu}=\inf\ee|X-Y|
$, where the infimum is is taken over all joint distributions of the random variables $X$ and $Y$  with marginals $\mu$ and $\nu$ respectively. We frequently use the following duality form of definition\footnote{$\operatorname{Lip}(h)$ is the Lipschitz constant of $h$ such that $|h(x_2)-h(x_1)|\leq \operatorname{Lip}(h)|x_2-x_1|$.},
\begin{equation*}
    \ws{\mu-\nu}:=\sup\left\{\left|\int_\xx h(x)d\mu(x)-\int_\xx h(x)d\nu(x)\right|,\;h\in C(\xx),\operatorname{Lip}(h)\leq 1\right\}. 
\end{equation*}
The discrete case, $\ws{\cdot}^d$, is nothing but to change the integral to summation. Let $B_W=\{\mu\in\psp(\xx): \ws{\mu-\delta_0}\leq 1\}$. 
Given a set $\G\subseteq\psp(\xx)$, we denote $\|\mu\|_\G=\inf_{\nu\in\G}\ws{\mu-\nu}$  by the distance from $\mu$ to $\G$, and $
    \G+r\bw:=\{\mu:\;\|\mu\|_\G\leq r\}%=\cup_{\mu\in\G}(\mu+r\bw). 
$\footnote{This is valid by definition.} by the $r$-neighborhood of $\G$.
\end{deff}
\begin{rem}
Note that $\mathbb{B}_W$ is dual to $\mathbb{B}_1$. For any $\mu\in \mathbb{B}_W$, the associated random variable $X$ should satisfy $\ee|X|\leq 1$, and vice versa. 
\end{rem}
The following well-known result estimates the  Wasserstein distance between two Gaussians.  

\begin{prop}\label{prop: compare}
Let $\mu\sim\mathcal{N}(m_1,\Sigma_1)$ and $\nu\sim\mathcal{N}(m_2,\Sigma_2)$ be two Gaussian measures on $\R^n$. Then
\begin{equation}\label{E: ws}
   |m_1-m_2|\leq \ws{\mu-\nu}\leq \left(\|m_1-m_2\|_2^2+\|\Sigma_1^{1/2}-\Sigma_2^{1/2}\|_F^2\right)^{1/2},%\footnote{The R.H.S. of \eqref{E: ws} is the $2^{\text{nd}}$-Wasserstein distance for two Gaussians.}.
\end{equation}
where $\|\cdot\|_F$ is the Frobenius norm. %We use   as the matrix norm. 
\end{prop}

%\begin{rem}
%For $\mu\sim\mathcal{N}(m_1,\Sigma)$ and $\nu\sim\mathcal{N}(m_2,\Sigma)$ with the same covariance, $\ws{\mu_1-\mu_2}=|m_1-m_2|$. 
%\end{rem}
On finite state spaces, total variation and Wasserstein distances manifest equivalence \cite[Theorem 4]{gibbs2002choosing}. We only show the following side of inequality in favor of the later proofs. 
\begin{prop}\label{prop: dis}
For any $\mu,\nu$ on some discrete and finite space $Q$, we have 
\begin{equation}
    \ws{\mu-\nu}^d\leq \operatorname{diam}(Q)\cdot\tv{\mu-\nu}.
\end{equation}
\end{prop}
Before proceeding, we define the set of transition probabilities of $\M_i$ from any box $[x]\subseteq\R^n$ as $$\T_i([x])=\{\mathcal{T}(x,\cdot):\;\mathcal{T}\in\trans{\mathcal{T}}_i,\;x\in[x]\},\;i=1,2,$$
and use the following lemma to approximate $\T_1([x])$. 
\begin{lem}\label{lem: inclu}
Fix any $\eps>0$, any box $[x]\subseteq\R^n$. For all $\kappa>0$, there exists a finitely terminated algorithm to compute an over-approximation of the set of (Gaussian) transition probabilities from $[x]$, 
such that
$$\T_1([x])\subseteq \wh{\T_1([x])}\subseteq\T_1([x])+\kappa\mathbb{B}_W, $$
%$$\mathcal{T}([x],\cdot)\subseteq \wh{\mathcal{T}([x],\cdot)}\subseteq\mathcal{T}([x],\cdot)+\kappa\mathbb{B}_W, $$
where $\wh{\T_1([x])}$ is the computed over-approximation set of Gaussian measures.
\end{lem}
\begin{rem}
The proof is completed in Appendix \ref{sec: appb}.
%Note that $\{\mathcal{N}(m\pm \kappa,s^2): \;\mathcal{N}(m,s^2)\in\T([x])\}\subseteq \T([x])+\kappa\mathbb{B}$.
%The latter ones have more choices of diffusion than the formal ones. 
The lemma renders the inclusions with larger Wasserstein distance to ensure no missing information about the covariances.
\end{rem}
\begin{deff}
For $i=1,2$, we introduce the modified transition probabilities for $\M_i=(\xx,\trans{\mathcal{T}}_i,x_0,\Pi,L)$ based on \eqref{E: T}. For all $\mathcal{T}_i\in\trans{\mathcal{T}}_i$, let \begin{equation}\label{E: mod}
    \ttil_i(x,\Gamma)=\left\{\begin{array}{lr} \mathcal{T}_i(x,\Gamma),\; \forall\Gamma\subseteq\W,\;\forall x\in W,\\
\mathcal{T}_i(x,\W^c),\;\Gamma=\partial\W,\;\forall x\in W,\\
1,\;\Gamma=\partial\W,\;x\in\partial\W.
\end{array}\right.
\end{equation}
Correspondingly, let $\tilde{\trans{\mathcal{T}}}$ denote the collection. Likewise, we also use $\tilde{(\cdot)}$ to denote the induced quantities of any other types w.r.t. such a modification. % any after-modification notations.
\end{deff}
\begin{rem}
We introduce the concept only for analysis. The above modification does not affect the law of the stopped processes since we do not care about the `out-of-domain' transitions.  We use a weighted point mass to represent the measures at the boundary, and the mean should remain the same. It can be easily shown that the Wasserstein distance between any two measures in $\tilde{\trans{\mathcal{T}}}(x,\cdot)$ is upper bounded by that of the non-modified ones. 
\end{rem}
\begin{thm}
For any $0\leq \ep_1<\ep_2$, we set $\M_i=(\xx, \tilde{\trans{\mathcal{T}}}_i, x_0,\Pi,L)$, $i=1,2$, where $\M_1$ is perturbed by point masses with intensity $\ep_1$, and $\M_2$ is perturbed by general $L_1$-perturbation with intensity $\ep_2$.
Then, under Assumption \ref{ass: part}, there exists a rectangular partition $Q$ (state-level relation $\alpha\subseteq\xx \times Q$),% that respects the boundary of labelling function,
a measure-level relation $\gamma_\alpha$ 
and a collection of transition matrices $\trans{\tta}$,  such that the system $\I=(Q,\trans{\tta},q_0,\Pi,L)$ abstracts $\M_1$ and is abstracted by $\M_2$ by the following relation:
\begin{equation}
    \M_1\preceq _{\gamma_\alpha}\I,\;\;\I\preceq_{\gamma_\alpha^{-1}}  \M_2.
\end{equation}
\iffalse
\begin{equation}
    \tilde{\trans{\mathcal{T}}}_1\preceq _{\gamma_\alpha}\trans{\tta},\;\;\trans{\tta}\preceq_{\gamma_\alpha^{-1}}\tilde{\trans{\mathcal{T}}}_2.
\end{equation}\fi
\end{thm}
\begin{pf}
We
    construct a finite-state IMC with partition $Q$ and an inclusion of transition matrices $\trans{\tta}$ as follows. By Assumption \ref{ass: part}, we use uniform rectangular partition on $\W$ and set $\alpha=\{(x,q): q=\eta\lfloor\frac{x}{\eta}\rfloor\}\cup\{(\Delta,\Delta)\}$, where $\lfloor\cdot\rfloor$ is the floor function and $\eta$ is to be chosen later. Denote the number of discrete nodes by $N+1$.
    
     Note that any family of (modified) Gaussian measures $\tmo_1$ is induced from  $\trans{\mathcal{T}}_1$ and should contain its information. %The $\gamma_\alpha$ is computed by the following steps. 
     For any $\tdt\in\tmo_1$ and $q\in Q$,
    \begin{enumerate}
        \item[(i)] %for any , %and for all $x\in\alpha^{-1}(q_i) $,  
        for all $\tnu\sim\tilde{\mathcal{N}}(m,s^2)\in\tilde{\T}_1(\alpha^{-1}(q),\cdot)$, store $\{(m_{l},\Sigma_{l})=(\eta\lfloor\frac{m}{\eta}\rfloor,\eta^2\fl{\frac{s^2}{\eta^2}})\}_l$;
        \item[(ii)] for each $l$, define $\tnu_{l}^\rff\sim\tilde{\mathcal{N}}(m_{l},\Sigma_{l})$ (implicitly, we need to compute $\nu_{l}^\rff(\Delta)$); compute $\tnu_{l}^\rff(\alpha^{-1}(q_j))$ for each $q_j\in Q\setminus\Delta$; %compute $\nu_{k}^\rff(\Delta)$; 
        \item[(iii)] for each  $l$,  define $\mu^{\text{ref}}_{l}=[\tnu_{l}^\rff(\alpha^{-1}(q_1)),\cdots,\tnu_l^\rff(\alpha^{-1}(q_N)),\tnu_l^\rff(\Delta)]$;
        \item[(iv)] compute
       $ \textbf{ws}:=(\sqrt{2N}+2)\eta 
        $ and $\textbf{tv}:=N\eta\cdot\textbf{ws} $;
        \iffalse
        \begin{equation}
            \operatorname{WS}=\frac{1}{2}\sqrt{\operatorname{tr}(\Sigma_i^{-1}(\Sigma_i+\eta)-I)+(n+1)^2\eta^2\sum\Sigma_i^{-1}-\log\operatorname{det}((\Sigma_i+\eta)\Sigma_i^{-1})}
        \end{equation}\fi
        \item[(v)] construct $\trans{\mu}=\bigcup_l\{\mu:\tv{\mu-\mu^{\text{ref}}_{l}}\leq \textbf{tv}(\eta),\; \mu(\Delta)+\sum_j^N\mu(q_j)=1 \}$;
        \item[(vi)]
        Let $\gamma_\alpha=\{(\tnu,\mu),\;\mu\in\trans{\mu}\}$ be a relation between $\tnu\in\tilde{\T}(\alpha^{-1}(q))$ and the generated $\trans{\mu}$.
        
        %Let row $i$ of the set-valued matrix $\trans{\tta}_\mathcal{T}$ be $\trans{\mu_i}$ for $i=1,\cdots,n$, and the $n+1$ row be $(0,\cdots, 1)$. 
    \end{enumerate}
    %Then,        $\gamma_{\iota}=\{(\tdt,\;\tta):\;\tta_i\in \iota_\alpha(\mathcal{T}(\alpha^{-1}(q_i),\cdot), \forall i\leq n\;\text{and}\; \tta_{n+1}=(0,\cdots,1)\}$ is the relation between transitions.
    Repeat the above step for all $q$, the relation $\gamma_\alpha$ is obtained.  The rest of the proof falls in the following steps. For $i\leq N$,  we simply denote $\G_i:=\tilde{\T}_1(\alpha^{-1}(q_i))$ and $\hat{\G}_i:=\wh{\tilde{\T_1}(\alpha^{-1}(q_i))}$.\\
\noindent
    \text{Claim 1}: For $i\leq N$, let $\trans{\tta_i}= \gamma_\alpha(\hat{\G}_i)$. Then the finite-state IMC $\I$ with transition  collection $\trans{\tta}$ abstracts $\M_1$.%I.e., we  construct $\tta_i$'s  for every $\tdt\in\tmo $ by $\iota_\alpha(\wh{\tdt([\alpha^{-1}(q_i)],\cdot)})$, then the collection $\trans{\tta}$ abstracts $\tmo$. %Note that $\iota_\alpha(\wh{\mathcal{T}(\cdot,\cdot)})$ and $\iota_\alpha(\wh{\tilde{\mathcal{T}}(\cdot,\cdot)})$ do not provide anything different.
    
    Indeed, for each $i=1,\cdots, N$ and each $\tdt$,  we have $\gamma_\alpha(\G_i)\subseteq \gamma_\alpha(\wh{\G}_i)$.  We pick any modified Gaussian $\tilde{\nu}\in\hat{\G}_i$, there  exists a $\tilde{\nu}^{\rff}$ such that (by Proposition \ref{prop: compare}) $\ws{\tilde{\nu}-\tilde{\nu}^{\rff}}\leq\ws{\nu-\nu^\rff}\leq \sqrt{2N}\eta$. 
    \iffalse
    \begin{equation*}
        \begin{split}
            \ws{\tilde{\nu}-\tilde{\nu}^{\rff}}\leq\ws{\nu-\nu^\rff}\leq \sqrt{2}\eta.
        \end{split}
    \end{equation*}\fi
    %By the above construction, any discrete $\mu\in\iota(\hat{G})$ is generated by comparison with some $\mu^\rff$, which is generated by some corresponding $\tilde{\nu}^{\rff}\in\hat{G}$ (see step (iii)).  
    We aim to find all discrete measures $\mu$ induced from $\tnu$ (such that their probabilities match on discrete nodes as requirement by Definition \ref{def: abs}).  %$\mu\in\gamma_\alpha(\hat{\G}_i)\subseteq\psp(Q)$ are related by . 
    All such $\mu$ should satisfy\footnote{Note that we also have $\ws{\mu-\mu^\rff}^d\leq  \ws{\mu-\tnu}^d+\ws{\tnu-\tnu^\rff}^d+\ws{\tnu^\rff-\mu^\rff}^d=\ws{\tnu-\tnu^\rff}^d$, but it is hard to connect $\ws{\tnu-\tnu^\rff}^d $ with $\ws{\tnu-\tnu^\rff}$ for general measures. This connection can be done if we only compare Dirac or discrete measures.},
    \begin{equation}
        \begin{split}
         \ws{\mu-\mu^\rff}^d & =   \ws{\mu-\mu^\rff}\\
         & \leq \ws{\mu-\tnu}+\ws{\tnu-\tnu^\rff}+\ws{\tnu^\rff-\mu^\rff}\\
                  &\leq (2+\sqrt{2N})\eta,
        \end{split}
    \end{equation}
   % $$\ws{\mu-\mu^\rff}^d =\sup_{h\in C(Q),\operatorname{Lip}\leq 1}\sum_{j=1}^n |h(q_j)[\tnu(\alpha^{-1}(q_j))-\tnu^\rff(\alpha^{-1}(q_j)]|\leq \ws{\tilde{\nu}-\tilde{\nu}^{\rff}}$$ %$ \ws{\mu-\mu^\rff} \leq \ws{\mu-\tilde{\nu}}+\ws{\tilde{\nu}-\tilde{\nu}^{\rff}}+\ws{\tilde{\nu}^{\rff}-\mu^\rff}$,
   where the first term of line 2 is bounded by, 
\begin{equation}\label{E: bound}
        \begin{split}
            \ws{\mu-\tilde{\nu}}&=\sup_{h\in C(\xx),\operatorname{Lip}(h)\leq 1}\left|\int_\xx h(x)d\mu(x)-\int_\xx h(x)d\tilde{\nu}(x)\right|\\
            &\leq \sup_{h\in C(\xx),\operatorname{Lip}(h)\leq 1}\sum_{j=1}^n\int_{\alpha^{-1}(q_j)}|h(x)-h(q_j)|d\tilde{\nu}(x)\\&\leq \eta \sum_{j=1}^n\int_{\alpha^{-1}(q_j)}d\tilde{\nu}(x)\leq \eta,
        \end{split}
    \end{equation}
    and the third term of line 2 is bounded in a similar way. 
    By step (v)(vi) and Proposition \ref{prop: dis}, 
     all possible  discrete measures $\mu$ induced from $\tnu$  should be included in  $\gamma_\alpha(\hat{\G}_i)$. Combining the above, for any $\tilde{\nu}\in\G_i$ and hence in $\hat{\G}_i$, there exists a discrete measures in $\tta_i\in\gamma_\alpha(\hat{\G}_i)$  such that for all $q_j$ we have  $\tilde{\nu}(\alpha^{-1}(q_j))=\tta_{ij}$. This satisfies the definition of abstraction.\\
    
    \noindent\text{Claim 2}:  $             \gamma_\alpha^{-1}(\gamma_\alpha(\G_i))\subseteq\G_i+(2\eta+N\eta\cdot\textbf{tv}(\eta))\cdot\bw
            $. This is to recover all possible (modified) measures $\tnu$ from the constructed $\gamma_\alpha(\G_i)$, such that their discrete probabilities coincide. 
            Note that, the `ref' information is recorded when computing $\gamma_\alpha(\G_i)$ in the inner parentheses. Therefore, for any $\mu\in\gamma_\alpha(\G_i)$ there exists a $\mu^\rff$ within a total variation radius $\textbf{tv}(\eta)$. We aim to find corresponding  measure $\tnu$ %\in\gamma_\alpha^{-1}(\gamma_\alpha(\G_i))\subseteq\psp(\xx)$ 
            that matches $\mu$ by their probabilities on discrete nodes. 
            All such $\tnu$ should satisfy,
            \begin{equation}\label{E: incl}
                \begin{split}
                    \ws{\tnu-\tnu^\rff}& \leq \ws{\tnu-\mu}+\ws{\mu-\mu^\rff}^d+\ws{\mu^\rff-\tnu^\rff}\\
                  &\leq 2\eta+N\eta\cdot\textbf{tv}(\eta),
                \end{split}
            \end{equation}
    where the bounds for the first and third terms are obtained in the same way as \eqref{E: bound}.
The second term is again by a rough comparison in Proposition \ref{prop: dis}. Note that $\tnu^\rff$ is already recorded in $\G_i$. The inequality in \eqref{E: incl} provides an upper bound of Wasserstein deviation between any possible satisfactory  measure and some $\tnu^\rff\in\G_i$. \\

\noindent\text{Claim 3}: If 
we can choose $\eta$ and $\kappa$ sufficiently small such that
\begin{equation}
    2\eta+N\eta\cdot\textbf{tv}(\eta)+\kappa\leq \ep_2-\ep_1,
\end{equation}
then $\I\preceq_{\gamma_\alpha^{-1}}  \M_2$.

Indeed, the $\trans{\tta}$ is obtained by $\gamma_\alpha(\hat{\G}_i)$ for each $i$. By Claim 2 and Lemma \ref{lem: inclu}, we have 
$$\gamma_\alpha^{-1}(\gamma_\alpha(\hat{\G}_i))\subseteq \hat{\G}_i+(2\eta+N\eta\cdot\textbf{tv}(\eta))\cdot\bw\subseteq \G_i+(2\eta+N\eta\cdot\textbf{tv}(\eta)+\kappa)\cdot\bw $$ for each $i$. %We denote $D_i^w:=\{\tilde{\mathcal{N}}(m\pm d,s^2):\;\tilde{\mathcal{N}}(m,s^2)\in\tilde{\T}(\alpha^{-1}(q_i)),d\leq w\} $ for each $i$. 
By the construction, we can verify that
$\tilde{\T}_2(\alpha^{-1}(q_i))=\G_i+(\ep_2-\ep_1)\cdot\bw$.  The selection of $\eta$ makes $\gamma_\alpha^{-1}(\gamma_\alpha(\hat{\G}_i))\subseteq\tilde{\T}_2(\alpha^{-1}(q_i))$, which completes the proof.
\end{pf}

\begin{rem}
The relation $\gamma_\alpha$ (resp. $\gamma_\alpha^{-1}$) provides a procedure to include all proper (continuous, discrete) measures that connect
with the discrete probabilities. The key point is to record $\tnu^\rff$, $\mu^\rff$, and the corresponding radius. These are nothing but finite coverings of the space of measures. This also explains the reason why we use `finite-state' rather than `finite' abstraction. The latter has a meaning of using finite numbers of representative measures to be the abstraction. 

To guarantee a sufficient inclusion,   conservative estimations, e.g. the bound $\sqrt{2N}\eta$ in Claim 1 and the bound in Proposition \ref{prop: dis}, are made. This estimation can be done  more accurately given more assumptions. For example, the deterministic systems (where $w$ becomes $\delta$) provide Dirac transition measures, the $\ws{\mu-\mu^\rff}^d=0$ and hence the second term in \eqref{E: incl} is $0$.
\end{rem}
\begin{rem}\label{rem: convex}
Note that, to guarantee the second abstraction based on $\gamma_\alpha^{-1}$, we search all possible measures that has the same discrete probabilities as $\mu\in\gamma_\alpha(\hat{\G}_i)$, not only those Gaussians with the same covariances as $\G_i$ (or $\hat{\G}_i$). Such a set of measures  provide a convex ball w.r.t. Wasserstein distance. This actually makes sense because in the forward step of creating $\I$, we have used both Wasserstein and total variation distance to find a convex inclusion of all Gaussian or Gaussian related measures. There ought to be some measures that are `non-recoverable' to Gaussians, unless we extract some `Gaussian recoverable' discrete measures in $\trans{\tta_i}$, but this loses the point of over-approximation. In this view, IMC abstractions provide unnecessarily larger inclusions than needed.

For the deterministic case, the above mentioned `extraction' is possible, since the transition measures do not have diffusion, the convex inclusion becomes a collection of  vertices themselves (also see Remark \ref{rem: tight}). Based on these vertices, we are able to use $\gamma_\alpha$ to find the $\delta$ measures within a convex ball w.r.t. Wasserstein distance. 

In contrast to the above special case \cite{liu2017robust}, where the uncertainties are bounded w.r.t. the infinity norm,  we can only guarantee the approximated completeness via a robust $\ls_1$-bounded perturbation with strictly larger intensity than the original point-mass perturbation. However, this indeed describes a general type of uncertainties for the stochastic systems to guarantee $\ls_1$-related properties, including probabilistic properties. Unless higher-moment specifications are of interests, uncertain $\ls_1$-random variables are what we need to be the analogue of perturbations in \cite{liu2017robust}.  
\end{rem}
\begin{cor}\label{cor: comp_inclu}
Given an LTL formula $\Psi$, let $S_i^{\nu_0}=\{\pp_X^{\nu_0}(X\vDash\Psi)\}_{X\in\M_i}$ ($i=1,2$) and $S_\I^{q_0}=\{\pim_I^{q_0}(I\vDash\Psi)\}_{I\in\I}$, where the initial conditions are such that $\nu_0(\alpha^{-1}(q_0))=1$. Then all the above sets are compact and $S_1^{\nu_0}\subseteq S_\I^{q_0}\subseteq S_2^{\nu_0}$.
\end{cor}
The proof in shown in Appendix \ref{sec: appb}.

\section{Conclusion}\label{sec: conclusion}
In this paper, we constructed an IMC abstraction for continuous-state stochastic systems with possibly bounded point-mass (Dirac) perturbations. We  showed that such 
abstractions are not only sound, in  the sense that the set of  satisfaction probability 
 of linear-time properties contains that of the original system, but
also approximately complete in the sense that the constructed IMC can be abstracted by another system with stronger but more general $\ls_1$-bounded 
perturbations. Consequently,
the winning set of the probabilistic specifications for a more perturbed continuous-state stochastic system contains 
that of the less Dirac perturbed system. Similar to most of the existing converse theorems, e.g. converse Lyapunov functions, the purpose is not to provide an efficient approach for finding
them, but rather to characterize the theoretical possibilities of having such existence.

It is interesting to compare with robust deterministic systems, where no random variables are involved. In \cite{liu2017robust}, both perturbed systems are  w.r.t. bounded point masses. More heavily perturbed systems abstract less perturbed ones and hence  preserve robust satisfaction of linear-time properties. However, when we try to obtain the approximated completeness via uncertainties in stochastic system, the uncertainties  should be modelled by more general $\ls_1$ random variables. Note that the probabilistic properties of random variables is dual to the weak topology of measures,  we study the measures and hence probability laws of processes instead of the state space \textit{per se}. The state-space topology is not sufficient to quantify the regularity of IMC abstractions. In contrast, the $\ls_1$ uncertain random variables is a perfect analogue of the uncertain point masses (in $|\cdot|$) for deterministic systems. If we insist on using point masses as the only type of uncertainties for stochastic systems, the IMC type abstractions would possibly fail to guarantee the completeness. For example, suppose the point-mass perturbations represents less precision of deterministic control inputs \cite[Definition 2.3]{majumdar2020symbolic}, the winning set decided by the $\ep_2$-precision stationary policies is not enough to cover that of the IMC  abstraction, which fails to ensure an approximated bi-similarity of IMCs compared to \cite{liu2017robust}.

For future work, it would be useful to extend the current
approach to robust stochastic control systems. It would be interesting to design algorithms to construct IMC (resp. bounded-parameter Markov decision processes) abstractions for more general robust stochastic (resp. control) systems with $\ls_1$ perturbations based on metrizable space of measures and weak topology. The size of state discretization can be refined given more specific assumptions on  system dynamics and linear-time objectives. For verification or control synthesis w.r.t. probabilistic safety or reachability problems, comparisons can be made with stochastic Lyapunov-barrier function approaches.

\newpage
\bibliographystyle{splncs04}
\bibliography{CAV}
\newpage
\appendix
\section{Proofs of Section \ref{sec: IMC}}\label{sec: appa}
\textbf{Proof of Corollary \ref{cor: marg}}.
\begin{pf}
It is clear that $Q$ under discrete metric is complete and separable. In addition, for each $t$, the space $(\Mm_t, \tv{\;\cdot\;})$ is complete and separable. By Lemma \ref{lem: compact}, each $(\Mm_t, \tv{\;\cdot\;})$ is also compact. For any sequence $\{\mu_n\}\subseteq\Mm_t$, a quick application of Theorem \ref{thm: proh} leads to the existence of a weakly convergent subsequence $\{\mu_{n_k}\}$ and a weak limit point $\mu$ in $\Mm_t$. By the definition of weak convergence and the discrete structure of $Q$, it is clear that for each $h\in C_b(\xx)$ and $t\in\mathbb{N}$, we have
$$\sum_\xx h(x)\mu_{n_k}(x)\rightarrow \sum_\xx h(x)\mu $$
in a strong sense, which concludes the compactness of $H$. Now we choose $\mu_1,\mu_2\in\Mm_t$, then $\alpha\mu_1+(1-\alpha)\mu_2\in\Mm_t$ for all $\alpha\in [0,1]$. Therefore, 
$$\alpha\sum_\xx h(x)\mu_1(x)+(1-\alpha)\sum_\xx h(x)\mu_2(x)=\sum_\xx h(x)[\alpha\mu_1+(1-\alpha)\mu_2](x)\in H$$
for all $\alpha\in [0,1]$. This shows the convex structure of $H$.
\end{pf}\\

\noindent\textbf{Proof of Proposition \ref{prop: law_weak}}.
\begin{pf}
We make a bit abuse of notation and 
define $\pi_T:  Q^{ \infty}\rightarrow \prod_0^T Q$ as the projection onto the finite product space of $Q$ up to time $T$. Since we do no emphasize the initial conditions, we also use $\pim$, $\Mm$ and $\Mm_t$ for short. By Tychonoff
theorem, any product of $Q$ is also compact w.r.t. the product topology. Therefore, any family of measures on $Q^T$ is tight and hence compact.  By Remark \ref{rem: prod}, for every $\pim\in\Mm$, we have
$\pim\circ \pi_T^{-1}=\otimes_{t=0}^T\mu_t$ (recall Remark \ref{rem: prod}) for some $\mu_t\in\Mm_t$, and $\{\pim\circ \pi_T^{-1}\}_{I\in\I}$ forms a compact set. Hence, every sequence $\{\pim_n\circ \pi_T^{-1}\}_n\subseteq \{\pim\circ \pi_T^{-1}\}_{I\in\I}$ with any finite $T$ contains a weakly convergent subsequence. We construct the convergent subsequence of $\{\pim_n\}_n$ in the following way.

We initialize the procedure by setting $T=0$. Then   $\Mm_0$ is compact, and there exists a weakly convergent subsequence $\{\pim_{0,n}\circ \pi_0^{-1}\}$. Based on $\{\pim_{0,n}\}$, we are able to see that it %$\{\pim_{0,n}\circ \pi_1^{-1}\}$ 
contains a weakly convergent subsequence, denoted by $\{\pim_{1,n}\}$, such that $\{\pim_{1,n}\circ \pi_1^{-1}\}$ weakly converges. By induction, we have $\{\pim_{k+1,n}\}\subseteq\{\pim_{k,n}\}$ for each $k\in \mathbb{N}$. Repeating this argument and picking the diagonal subsequence $\{\pim_{n,n}\}$, then $\{\pim_{n,n}\}$ has the property that  $\{\pim_{n,n}\circ \pi_T^{-1}\}$ is weakly convergent for each $T$ . We denote the weak limit point of each $\{\pim_{n,n}\circ \pi_T^{-1}\}$ by $\otimes_{t=0}^T\mu_t$. By the way of construction, we have 
$$\otimes_{t=0}^T\mu_t(\cdot)=\otimes_{t=0}^{T+1}\mu_t(\cdot\times Q),\;\;\forall T\in\mathbb{N}.$$ 
By Kolmogorov's extension theorem, there exists a unique $\pim$ on $Q^\infty$ such that $\otimes_{t=0}^T\mu_t(\cdot)=\pim\circ\pi_T^{-1}(\cdot)$ for each $T$. 

We have seen that for each $\{\pim_n\}$, the constructed subsequence satisfies $\pim_{n,n}\Rightarrow\pim$, which concludes the claim. 
\end{pf}\\

\noindent\textbf{Proof of Theorem \ref{thm: imc_bound}}
\begin{pf}
Since we do not emphasize the initial conditions, we simply drop the superscripts $q_0$ for short.
Given $I\in\I$ with any initial condition, %the event set $\ob=Q^\infty$, and $\iim$,
the corresponding canonical space is $(\ob,\ffb,\pim_I)$. By Proposition \ref{prop: law_weak}, every sequence $\{\pim_n\}\subseteq \Mm$ has a weakly convergent subsequence, denoted by $\{\pim_{n_k}\}$, to a $\pim\in\Mm$ of some $I$. Note that for any $I$, the measurable set $\{I\vDash\Psi\}=\{\w: \w\vDash\Psi\}\in\ffb$ is the same due to the identical labelling function. It is important to notice that due to the discrete topology of $\ob$, every Borel measurable set $A\in\ffb$ is such that $\partial A=\emptyset$.  By Definition \ref{def: weak_conv} we have $\pim_{n_k}(I_{n_k}\vDash\Psi)\rightarrow \pim(I\vDash\Psi)$ for all $\Psi$. The compactness of $S^{q_0}$ follows immediately. The convexity of the set of laws is based on the tensor product of convex polytopes \cite{valby2006category}. To show the convexity of $S^{q_0}$, we notice that, for any $q_0,\cdots q_{n_t}\in Q$ and $I\in\I$,
\begin{equation*}
    \begin{split}
       &\pim_I\left(I_0=q_0, \cdots, I_t=q_{n_t}, I_{t+1}=q_{n_{t+1}}\right)\\
       \in&\{\tta_{n_{t+1},n_t}\tta_{n_t,n_{t-1}}\cdots\tta_{n_1,0}\delta_{q_0}: \tta\in [\![\tta]\!]\}  
    \end{split}
\end{equation*}
and hence forms a convex set. Immediately, the convexity holds for\\ $\{\pim_I(\prod_{i=1}^n \Gamma_i)\}_{I\in\I}$
for any cylinder set $\prod_{i=1}^n \Gamma_i $. By a standard monotone class argument, $\{\pim_I(A)\}_{I\in\I}$ is also convex for any Borel measurable set $A\in \ffb$, which implies the convexity of $S^{q_0}$ in the statement. %Note that for discrete topology, every $h:Q^\infty\rightarrow \mathbb{Z}$ is continuous. Therefore, indicator functions are continuous and
%$\sum_\Omega \mathds{1}_{\{\iim\vDash\Psi\}}\pim^\iim(\w)=\pim^\iim(\iim\vDash\Psi),\;\;\forall \iim.$
\end{pf}\\

\noindent\textbf{Proof of Proposition \ref{cor: compact}}
\begin{pf}
Note that the laws are associated with $X$ with $X_0=x_0$, which actually means the stopped process $X^\tau$ (Recall notations in Section \ref{sec: sde}). Now that $X_{t\wedge\tau}\in \overline{\W}$ for each $t$, the state space of $X$ is compact, so is the countably infinite product. By a similar argument as Proposition \ref{prop: law_weak}, we can conclude the first part of the statement. Note that by assumption, the partition respects the boundary of the labelling function. Hence, for all formula $\Psi$, the boundary of $\{X\vDash\Psi\}\in\mathbf{F}$ has measure $0$. The second part can be concluded directly by Definition \ref{def: weak_conv}.
\end{pf}\\

\noindent\textbf{Proof of Lemma \ref{lem: connect}}
\begin{pf}
Note that $X$ is on $(\Omega,\ff,\pp_X^{\nu_0})$ and $I$ is on $(\ob,\ffb,\pim_I^{q_0})$. We first show the case when $\nu_0=\delta_{x_0}$ for any $x_0\in q_0$. That is, for  $X_0=x_0$ a.s. with any $x_0\in q_0$, there exists a unique law of some $I\in\I$  such that  $\pp_X^{x_0}(X\vDash\Psi)=\pim_I^{q_0}(I\vDash\Psi)$ for any $\Psi$.

Let $\nu_t$ denote the marginal distribution of $\pp_X^{x_0}$ at each $t$. Let $\Mm_t=\{\mu_t\}_{I\in\I}$ denote the set of marginal distributions of $\{\pim_I^{q_0}\}_{I\in\I}$. %It is enough to just look at transitions from $x\in\W$. 
Now, 
at $t=1$,  $\nu_1(q_j)=\mathcal{T}(x_0,q_j)\delta_{x_0}$ for all $j\in\{1,2,\cdots, N+1\}$. Suppose $q_0$ is the $i^{th}$ element of $Q$, by the construction of IMC, we have 
$$ \undert_{ij} \leq\nu_1(q_j)=\int_{q_j}\delta_{x_0}\mathcal{T}(x_0,dy) \leq \overt_{ij} ,\;\;\forall x_0\in q_0\;\text{and}\;\forall j\in\{1,2,\cdots, N+1\}. $$
Since $\sum_{q\in Q} \nu_1(q)=1$, by letting $\mu_1=(\nu_1(q_1),\nu_1(q_2),\cdots,\nu_1(q_{N+1}))^T$, we have automatically $\mu_1\in\Mm_1$ by definition. Note that $\mu_1$ is unique w.r.t. $\tv{\;\cdot\;}$, and has the property that $\mu(q)=\nu(q)$ for each $q\in Q$.

Similarly, at $t=2$, we have 
$$\undert_{ij}\mu_1(q_i)\leq \int_{q_j}\int_{q_i}\nu_1(dx)\mathcal{T}(x,dy)\leq \overt_{ij}\mu_1(q_i),\;\forall i,j\in\{1,2,\cdots, N+1\}, $$
where $\mathcal{T}$ may not be the same as that of $t=1$. Therefore, for any $j\in \{1,2,\cdots,N+1\}$, $$\nu_2(q_j)=\sum_{i=1}^{N+1}\int_{q_j}\int_{q_i}\nu_1(dx)\mathcal{T}(x,dy)$$
and there exists a $\mu_2$ such that $\sum_i\undert_{ij}\mu_1(q_i)\leq \mu_2(q_j)=\nu_2(q_j)\leq \sum_i\overt_{ij}\mu_1(q_i)$, which means (by \eqref{E: meas}) $\mu_2\in\Mm_2$. In addition, there also exists a $\pim^{q_0}$ such that its one-dimensional marginals up to $t=2$ admit $\mu_1$ and $\mu_2$, and
satisfies
$$\pim^{q_0}[I_0=q_0,I_1=q_i,I_2=q_j]=\pp_X^{x_0}[X_0=x_0,X_1\in q_i,X_2\in q_j]. $$
Repeating this procedure, there exists a unique $\mu_t\in\Mm_t$ w.r.t. $\tv{\;\cdot\;}$ for each $t$, such that $\mu_t(q)=\nu_t(q)$ for each $q\in Q$. It is also clear that for each given $x_0\in q_0$ and each $t$, the selected $\pim^{q_0}$ satisfies%$\mu_t$ 
$$\pp_X^{x_0}(\prod_0^{t-1} A_i)%=\otimes_0^{t-1}\nu_i(\prod_0^{t-1} A_i)
=\pim^{q_0}(\prod_0^{t-1} A_i)=\pim^{q_0}(\prod_0^{t-1} A_i\times Q) ,\;\;A_i\in\mathscr{B}(Q).$$
By Kolmogrov extension theorem, there exists a unique law $\pim_I^{q_0}$ of some $I\in\I$ such that each $T$-dimensional distribution coincides with $\otimes_0^T\mu_i$, and, for each given $x_0\in q_0$, 
$\pp_X^{x_0}(\Gamma)=\pim_I^{q_0}(\Gamma)$ for all $ \Gamma\in \mathscr{B}(Q^\infty)=\ffb.$
Due to the assumption that $L(x)=L_\I(q)$ for all $x\in q$ and $q\in Q$, we have 
$$\{L^{-1}(\Psi)\}=\{L_\I^{-1}(\Psi)\}\in \ffb$$%\subsetneq\ff$$
for all LTL formula $\Psi$, which implies $\{X\vDash\Psi\}=\{I\vDash\Psi\}$ by definition.
Thus, given  $x_0\in q_0$, the above  $\pim_I^{q_0}$   should satisfy
$\pp_X^{x_0}(X\vDash\Psi)=\pim_I^{q_0}(I\vDash\Psi)$.

Based on the above conclusion, as well as the definition of $\pp_X^{\nu_0}$ and the convexity of $S^{q_0}$ (recall Theorem \ref{thm: imc_bound}), the result for more general initial distribution $\nu_0$ with $\nu_0(q_0)=1$ can be obtained.
\end{pf}\\

\noindent\textbf{Proof of Proposition \ref{prop: converge}}.
\begin{pf}
Let $\mu,\nu\in\Mm_t$ for each $t$, and $V,K\in [\![\tta]\!]$ be any stochastic matrices generated by $\I$. Then, for each $t$, we have
\begin{equation}
    \begin{split}
        \tv{V^T\mu-K^T\nu}&\leq \tv{V^T\mu-V^T\nu}+\max_i\tv{V_i-K_i}\tv{\nu}\\
        &\leq \frac{1}{2}\max_{i,j}\tv{V_i-K_j}\tv{\mu-\nu}+\eps\tv{\nu}\\
        &\leq \tv{\mu-\nu}+\eps.
    \end{split}
\end{equation}
This implies that the total deviation of any $\tilde{\mu},\tilde{\nu}\in\Mm_{t+1}$ is bounded by $$\max_{\Mm_t}\tv{\mu-\nu}+\eps.$$
Note that at $t=0$, $\max_{\Mm_0}\tv{\mu-\nu}=0$. Hence, at each $t>0$, as $\eps\rightarrow 0$, $$\max_{\Mm_t}\tv{\mu-\nu}\rightarrow 0. $$
By the product topology and Kolmogrov extension theorem, the set $\{\pim\}_{I\in\I}$ is reduced to a singleton. The conclusion follows after this.
\end{pf}

\section{Proofs of Section \ref{sec: complete}}\label{sec: appb}
\textbf{Proof of Lemma \ref{lem: inclu}}.
\begin{pf}
It can be proved, for example, using inclusion functions. Let $\mathbb{IR}^n$ denote the set of all boxes in $\R^n$. Let $[f]:\mathbb{IR}^n\rightarrow\mathbb{IR}^n$ be a convergent inclusion function of $f$ satisfying (i) $f([x])\subseteq[f]([x])$ for all $[x]\in\mathbb{IR}^n$; (ii) $\lim_{\lambda([y])\rightarrow0 }\lambda([f]([x])=0 $, where $\lambda$ denote the width.
Similarly, 
let $[B]:\mathbb{IR}^n\rightarrow\mathbb{IR}^{n\times k}$ be a convergent inclusion matrix of $b(x)$ and satisfy (i) $b([x])\subseteq[B]([x])$ for all $[x]\in\mathbb{IR}^n$; 
(ii) $\lim_{\lambda([B])\rightarrow0 }\lambda([B]([x])=0 $, where $\lambda([B]([x]):=\max_{i,j}\lambda([B_{ij}])$.

Without loss of generality, we assume that $\kappa<1$. Due to the Lipschitz continuity of $f$ and $b$, we can find inclusions such that $\lambda([f]([y])\leq L_f\lambda([y])$ for any subintervals of $[x]$, and similarly $\lambda([B]([y])\leq L_b\lambda([y])$. For each such interval $[y]$, we can obtain the interval 
%$$[z]=\left[([f]([y])+\ep\mathbb{B})^2+[B]([y])[B]^*([y])\right]^{1/2}. $$
$[m]=[f]([y])$
and 
$[s^2]=[B]([y])[B]^*([y])$.
Let $T$ denote the collection of Gaussian measures with mean and covariance of all such intervals ($[m]$ and $[s^2]$), and $\wh{\T_1([x])}$ be its union. Then $\wh{\T_1([x])}$ satisfies the requirement. Indeed, we have $\T_1([x])\subseteq \wh{\T_1([x])}$. For the second part of inclusion, we have for any $\mu\in \wh{\T_1([x])}$ and $\nu\in \T_1([x])$,
\begin{equation}
   \ws{\mu-\nu}^2\leq [m]^2+[s^2]\leq (L_f\lambda([y]))^2+(L_b\lambda([y]))^2,
\end{equation}
where we are able to choose $[y]$ arbitrarily small such that $(L_f\lambda([y]))^2+(L_b\lambda([y]))^2<\kappa^2$. The second part of inclusion can be completed by such a choice of $[y]$.
\end{pf}\\

\noindent\textbf{Proof of Corollary \ref{cor: comp_inclu} }
\begin{pf}
The first part of the proof is provided in Section \ref{sec: IMC}. The second inclusion is done in a similar way as Lemma \ref{lem: connect} and Theorem \ref{thm: inclusion}. Indeed, by the definition of abstraction, for any $\mu\in\Mm_t$, there exists a marginal measure $\nu$ of some $X\in\M_2$ such that their probabilities match on discrete nodes. By the same induction as Lemma \ref{lem: connect}, we have that for any law $\pim_I^{q_0}$ of some $I\in\I$, there exists a $\pp_X^{\nu_0}$ of some $X\in\M_2$ such that the probabilities of any $\Gamma\in\mathscr{B}(Q^\infty)$ match. The second inclusion follows after this. The compactness also follows a similar way as Proposition \ref{cor: compact}. Note that, $S_1$ may not be convex, but $S_\I$ and $S_2$ are (also see details in Remark \ref{rem: convex}).
\end{pf}
%
\iffalse

\fi
\end{document}